\documentclass[reqno,12pt]{amsart}
\usepackage{amsmath, amssymb, amsthm,amsfonts} 
\usepackage[english]{babel}
\usepackage{bbm,bm}
\usepackage{graphicx}
\usepackage{url}
\usepackage{epstopdf}
\usepackage[a4paper,bindingoffset=0.2cm,left=1cm,right=1cm,top=2.5cm,bottom=2cm,footskip=.8cm]{geometry}
\usepackage{subcaption}
\usepackage{booktabs}
\usepackage{textcomp}
\usepackage{longtable}
\usepackage{pdflscape}
\usepackage[algo2e,ruled,ruled,vlined,linesnumbered]{algorithm2e}
\usepackage{multicol}
\usepackage{multirow}
\usepackage{wrapfig}
\usepackage{diagbox}
\usepackage{floatrow}
\usepackage[dvipsnames]{xcolor}
\usepackage[]{algorithm2e}
\usepackage{lscape}
\usepackage{mathrsfs} 
\usepackage{rotating}
\usepackage{amsbsy,enumerate}
\usepackage{graphicx}
\usepackage{comment}
\usepackage{mathrsfs} 
\usepackage{xcolor}
\usepackage{tikz}
\usepackage{pgfplots}
\usepackage{arydshln}

\newcommand{\ee}{\mathbb{E}}
\newcommand{\pp}{\mathbb{P}}

\newcommand{\vb}{\vspace{4mm}}

\newcommand{\scv}{\textsc{scv}}

\newcommand{\objAS}{{\mathscr L} ( {\boldsymbol x})}

\newcommand{\Rene}[1]{\textcolor{blue}{[Rene: #1]}} 

\allowdisplaybreaks

\title[Accurate and efficient optimization of large-scale appointment schedules]{Accurate and efficient approximation of\\ large-scale appointment schedules}
\author{R. Bekker, B. Bharti, and M. Mandjes}

\begin{document}

\begin{abstract}

Setting up optimal appointment schedules requires the computation of an inherently involved objective function, typically requiring distributional knowledge of the clients' waiting times and the server's idle times (as a function of the appointment times of the individual clients). 
A frequently used idea is to approximate the clients' service times by their phase-type counterpart, thus leading to explicit expressions for the waiting-time and idle-time distributions. 
This method, however, requires the evaluation of the matrix exponential of potentially large matrices, which already becomes prohibitively slow from, say, 20 clients on.

\noindent 
In this paper we remedy this issue by recursively approximating the distributions involved relying on a two-moments fit. 
More specifically, we approximate the sojourn time of each of the clients by a low-dimensional phase-type, Weibull or Lognormal random variable with the desired mean and variance. 
Our computational experiments show that this elementary, yet highly accurate, technique facilitates the evaluation of optimal appointment schedules even if the number of clients is large. 
The three ways to approximate the sojourn-time distribution turn out to be roughly equally accurate, except in certain specific regimes, where the low-dimensional phase-type fit performs well across all instances considered.
As this low-dimensional phase-type fit is by far the fastest of the three alternatives, it is the approximation that we recommend.

\vb
\noindent
{\sc Keywords.} Appointment scheduling $\circ$ two-moments fit $\circ$ phase-type distribution $\circ$  queueing $\circ$ approximation

\vb

\noindent
{\sc Affiliations.} \emph{Bharti} and \emph{Michel Mandjes} are with the Korteweg-de Vries Institute for Mathematics, University of Amsterdam, Science Park 904, 1098 XH Amsterdam, the Netherlands. MM is also with Mathemetical Institute, Leiden University, Leiden, the Netherlands; E{\sc urandom}, Eindhoven University of Technology, Eindhoven, the Netherlands; and Amsterdam Business School, Faculty of Economics and Business, University of Amsterdam, Amsterdam, the Netherlands. 

\noindent The research of B and MM was supported by the European Union’s Horizon 2020 research and innovation programme under the Marie Skłodowska-Curie grant agreement no.\ 945045, and by the NWO Gravitation project NETWORKS under grant no.\ 024.002.003. \includegraphics[height=1em]{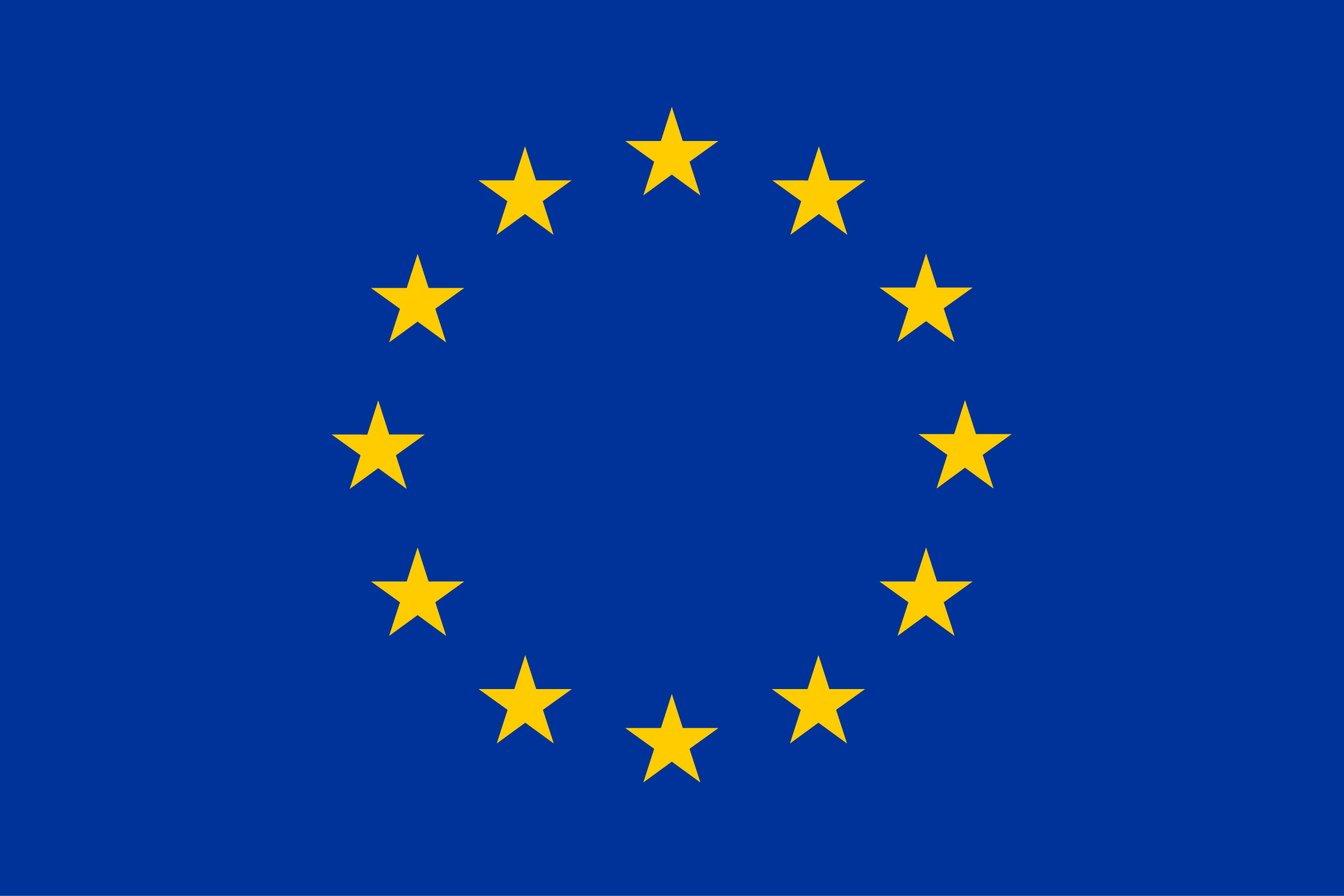}

\noindent
\emph{Ren\'e Bekker} is with the
Department of Mathematics,
Vrije Universiteit Amsterdam,
De Boelelaan 1111,
1081 HV Amsterdam,
the Netherlands. 

\vb

\noindent
\end{abstract}

\maketitle


\section{Introduction}\label{INTRO}
Appointment-based service systems arise in a broad variety of settings, such as in healthcare and parcel delivery  \cite{ahmadi2017, kuiper2017, kuiper2023, wang1997,zhan2021}. 
The main objective of appointment schedules is to find a sound balance between the interests of the clients and the service provider. 
This is achieved by employing an objective function that incorporates both clients' waiting times and the server's idle times.  
A complicating issue lies in the inherently random nature of the clients' service times. 

In mathematical terms, the objective function to be minimized can be defined as follows. Suppose there are $n$ clients to be scheduled, in a given order, and that we know the distribution of their service times. 
Let $I_j$ be the idle time {\it before} the arrival of the $j$-th client, and $W_j$ the corresponding waiting time. 
Define by ${\boldsymbol x}$ the vector of interarrival times, having dimension $n-1$. 
This concretely means that $x_j:= t_{j+1}-t_j\geqslant 0$ is the length of the time interval between the appointment time $t_j$ of client $j$ and the appointment time $t_{j+1}$ of client $j+1$. 
Our aim is to minimize the objective function, sometimes referred to as the {\it loss function},
\begin{equation}\label{OBJ0} {\mathscr L} ( {\boldsymbol x})  =  \omega \sum_{j=1}^n {\mathbb E}\,I_{j} + (1-\omega) \sum_{j=1}^n {\mathbb E}\,W_{j},\end{equation}
over the componentwise positive vector ${\boldsymbol x}$; observe that $I_j$ and $W_j$ depend on $x_1,\ldots,x_j$. 
The weight $\omega\in[0,1]$ reflects the relative importance of the idle times relative to that of the waiting times (and can be set by the service provider, based on her perception of this relative importance). 
A `generous' schedule, characterized by relatively long interarrival times, favors the clients because it leads to low waiting times. In contrast, a `tight' schedule, marked by relatively short interarrival times, favors the server because it leads to low idle times. 

In the literature the above optimization problem has been widely studied. 
An inherent problem is that, for the clients' service times having given distributions, typically the mean idle times and mean waiting times cannot be given by closed-form expressions. 
This led to the idea of approximating these service times by appropriately chosen {\it phase-type} distributions. These distributions have the potential to accurately approximate any distribution on the positive real numbers \cite[Thm.\ III.4.2]{asmussen2003}, while also allowing explicit calculations \cite{wang1997}. 
In particular, closed-form expressions for the mean idle and waiting times can be given, albeit in terms of matrix inverses and matrix exponentials; see e.g.\ \cite{kuiper2015}.

The phase-type approach has an intrinsic drawback as well: the dimension of the objects that need to be evaluated for the approximation of ${\mathscr L} ({\boldsymbol x})$ can become prohibitively large. 
More concretely, with $d_j$ denoting the dimension of the underlying phase-type approximation pertaining to the service time of the $j$-th client and with $D:=d_1+\cdots+d_n$, the method works with $D\times D$ matrices. 
Because of the matrix inverses and matrix exponentials needed to evaluate ${\mathscr L} ( {\boldsymbol x})$, this method effectively scales as $O(D^3)$. 
When minimizing ${\mathscr L} ( {\boldsymbol x})$ over ${\boldsymbol x}$, for instance by applying a standard Newton-type minimization routine, the objective function needs to be frequently evaluated (i.e., for different arguments ${\boldsymbol x}$), rendering this approach infeasible for large values of $n$ and/or $D$. 
In practice this means that this method cannot be used from roughly $n=20$ on, in particular in cases that some of the service times are of relatively low variability (corresponding to a phase-type distribution of high dimension). 
Thus, there is a clear need for a fast and accurate algorithm to determine optimal appointment schedules that also works when $n$ and/or $D$ are relatively large. 
The main objective of this paper is to develop such an algorithm. 

\vb

We now elaborate on the contributions of our paper. 
As has been observed in the literature \cite{kuiper2023}, both the mean idle times and mean waiting times can be evaluated if the {\it sojourn time} distribution of each individual client is known, with the sojourn time being defined as the sum of the waiting time and the service time. Indeed, with $R_i$ denoting the sojourn time of the $i$-th client, it requires some elementary arguments to verify that we can rewrite the loss function ${\mathscr L} ( {\boldsymbol x})$ to (with $y^+ = \max \{y,0\}$ for $y\in{\mathbb R}$)
\begin{equation}
    \label{OBJ}
{\mathscr L} ( {\boldsymbol x})  =  \omega \sum_{j=1}^n {\mathbb E}(x_{j-1}-R_{j-1})^+
+(1-\omega) \sum_{j=1}^n {\mathbb E}(R_{j-1}-x_{j-1})^+.\end{equation}
Our approach is to find a fast and accurate approximation scheme for the distributions of the sojourn times, as this facilitates the approximation of ${\mathscr L} ( {\boldsymbol x})$. 
These sojourn times obey a recursive relation, commonly known as the {\it Lindley recursion}:  for $i=0,1,\ldots$,
\begin{equation}\label{LR}
R_{i+1} = \max\{R_{i}-x_i,0\}+B_{i+1},    
\end{equation}
where $R_0=x_0=0.$ Let the random vector ${\boldsymbol B}=(B_1,\ldots,B_n)$ consist of independent non-negative components, representing the clients' service times. 
For evident reasons, the random variables $R_{i}$ and $B_{i+1}$ can be assumed independent. 
Throughout this paper we denote $\beta_i:={\mathbb E}\,B_i$ and $\sigma^2_i:={\mathbb V}{\rm ar}\,B_i$, for $i\in\{1,\ldots,n\}.$ 
Then our goal is to set up an efficient and accurate recursive procedure that approximates the mean and variance of the $R_i$ in terms of the mean service times $(\beta_1,\ldots,\beta_i)$ and the variances of the service times $(\sigma_1^2,\ldots,\sigma_i^2)$, for given interarrival times ${\boldsymbol x}.$
The underlying idea is that we approximate the sojourn times of each of the clients by  convenient random variables with the desired means and variances, i.e., a `two-moments fit'. We consider three classes of two-moments fits: one based on a (low dimensional) phase-type distribution, one based on the Weibull distribution, and one based on the Lognormal distribution. We specifically point out how the objective function ${\mathscr L}({\boldsymbol x})$ can be evaluated for each of the three fits. 

The main contribution of this work is a method by which we can evaluate the objective function with a computational effort that scales as just $O(n)$, rather than $O(D^3)$. 
We thus achieve a huge speedup, but remarkably, as indicated by our numerical experiments, the loss of accuracy is minimal. 
The experiments in addition show that the substantial reduction of the computation time allows the generation of optimal schedules even for rather large numbers of clients. 
Comparing the three two-moments fits, we conclude that the approximation based on a phase-type two-moments fit offers the best performance when taking both accuracy and execution time into account; the approaches based on the Weibull fit and the Lognormal fit are roughly equally accurate (except in specific regimes), but significantly slower.

\vb

We proceed by providing a brief account of the related literature. 
The approach that is probably closest to ours, is \cite{de1989}. 
While we present a more detailed comparison between the two approaches in Subsection~\ref{COMP_DEK}, we already mention that \cite{de1989} primarily focuses on evaluating steady-state quantities, while the metrics we target are inherently transient.
Algorithms of this type have a broad application potential, as witnessed by the sizeable number of papers in which the algorithm of \cite{de1989} has been used. 
It has in particular been widely applied in transport, logistics and inventory studies; see for instance \cite{klosterhalfen2018,wagner2004,zijm1994}.

The literature on approximations for single-server queues being extensive, we restrict ourselves to discussing a number of key references. 
For the more general class of GI/G/1 queues exact results on its transient performance are predominantly in terms of multi-dimensional transforms and characteristic functions. 
We refer, for instance, to the results in \cite[Chapters I and II]{prabhu1998}, where specifically \cite[Theorem I.29]{prabhu1998} presents various relevant multi-dimensional transforms. 
Although reliable numerical inversion algorithms are available \cite{abate1995}, their performance substantially degrades when it comes to double or triple inversion, rendering this approach essentially infeasible. 
A pragmatic remedy is to replace the interarrival times and service times by their phase-type counterparts, leading to explicit expressions for e.g.\ the queue-length and workload distributions, albeit in terms of possibly rather involved matrix operations; see the approaches followed in e.g.\ \cite{neuts1981,seelen1984}. 
It was also recognized that assuming that the interarrival times are of phase-type and the service times have a general distribution, or vice versa, facilitates a fairly clean analysis, but still in the transform domain; see e.g.\ \cite{bux1979}, and, in the terminology of the model's ruin-theoretic counterpart, \cite{lewis2008}. 
In the regime that the system load is above (say) 95\%, one often relies on the explicit heavy-traffic approximation; see e.g.\ \cite[Section X.7]{asmussen2003}. 
It is important to note that the vast majority of the above references exclusively consider steady-state performance measure (inherently assuming stochastically homogeneous customers), while the evaluation of our objective function 
${\mathscr L}({\boldsymbol x})$, as given in \eqref{OBJ}, requires transient analysis (allowing heterogeneous customers). 

In the specific context of appointment scheduling, exact evaluation methods of the objective function ${\mathscr L}({\boldsymbol x})$ are only known for phase-type service times \cite{kuiper2015,wang1997}. In case of other service times, one typically relies on simulation-based techniques \cite{ho1992,robinson2003}. Although simulation-based optimization is a well-developed branch of the literature \cite[Chapter VIII]{asmussen2007}, clearly, a method that minimizes ${\mathscr L}({\boldsymbol x})$ without the need to resort to simulation is to be preferred. 

\vb

This paper is organized as follows. In Section \ref{sec:loss} we outline the approach to approximate our objective function ${\mathscr L} ( {\boldsymbol x})$, whereas Section \ref{sec:alg} gives the specific algorithms we use to approximate it. We present extensive numerical experiments in Section \ref{NUM}, corroborating the claim that the proposed approach outperforms, by a huge margin, existing alternatives in terms of speeds, while giving up a minimal amount of accuracy. A discussion and concluding remarks are given in Section \ref{DIS}.

\section{Objective function and recursion} \label{sec:loss}


In this section, we give an outline of our approach to approximate the objective function $\objAS$, as was defined in Equation~(\ref{OBJ}). In addition, we point out how our methodology differs from the related approach followed in \cite{de1989}.

\subsection{Outline of our approach}
When aiming at the computation of \eqref{OBJ}, a useful key identity is as follows; for weights $A$ and $B$ and a constant $x\geqslant 0$, it evidently holds that
\begin{equation*}
    A(X-x)1\{X\geqslant x\}+B(x-X)1\{X<x\}=(A+B)(X-x)1\{X\geqslant x\}+B(x-X).
\end{equation*}
The loss function can therefore be written as 
\begin{equation*}
      {\mathscr L} ( {\boldsymbol x})=  \sum_{i=1}^n 
      {\mathscr L}_i ( {x}_{i-1}),\:\:\:\: \text{with} \:\:\:\:
      {\mathscr L}_i ( {x}):=
      {\mathbb E}\big(R_{i-1} - x\big)^+ + \omega\, \big(x - {\mathbb E} R_{i-1}\big).
\end{equation*}
This means that, to evaluate the objective function, we just need to compute $\ee(R_i - x)^+$; note that $\ee R_i$ then follows as a special case (i.e., by inserting $x=0$). 

The complication that we need to overcome is that no closed-form expressions for $\ee(R_i - x)^+$ (and $\ee R_i$) are available in the case of general service times $B_i$, so that we have to resort to approximations. 
The concept advocated in this paper, is to approximate $\ee(R_i - x)^+$ by
recursively determining the mean $r_i$ and variance $v_i$ of $R_i$, and using a `two-moment fit' for the sojourn time.
Applying the Lindley recursion~(\ref{LR}) and the fact that $R_i$ and $B_{i+1}$ are independent, the mean and variance satisfy the relations
\begin{equation}\label{Ri}
r_{i+1} := {\mathbb E}\,R_{i+1}
 = \ee(R_i - x_i)^+ + \beta_{i+1}
 = \int_{0}^{\infty} y\, f_{R_i}(y+x_i)\,{\rm d}y + \beta_{i+1},
\end{equation}
and
\begin{eqnarray}
v_{i+1} := {\mathbb V}{\rm ar}\,R_{i+1}
 &=& {\mathbb V}{\rm ar}\,(R_i - x_i)^+ + \sigma_{i+1}^2 \nonumber \\
 &=& \int_{0}^{\infty} y^2\,f_{R_i}(y+x_i)\,{\rm d}y-\left(\int_{0}^{\infty} y\,f_{R_i}(y+x_i)\,{\rm d}y\right)^2 +\sigma_{i+1}^2,    \label{Vi}
\end{eqnarray}
respectively, with $f_{R_i}(\cdot)$ the density of $R_i$, assuming it to exist. 
Hence, to determine $r_{i+1}$ and $v_{i+1}$, so as to be able to evaluate ${\mathscr L}_i ( {x})$, we essentially need to compute the quantities $\ee(R_i - x_i)^+$ and ${\mathbb V}{\rm ar}\,(R_i - x_i)^+$.
Our main idea is that we approximate the sojourn time $R_i$ by a conveniently chosen, relatively simple random variable.
More specifically, in our approximation we take the sojourn time from one of the following three types: a (i) low-dimensional phase type~(Ph), (ii) Weibull~(W), or (iii) Lognormal~(LN) distribution. 
A key advantage of these classes of distributions is that they allow for a tractable two-moment fit, whereas they additionally enable closed-form expressions of (\ref{Ri}) and (\ref{Vi}). 
In Subsection~\ref{subs:PH} we will further elaborate on the motivation for phase-type distributions. Note that in (\ref{LR}) the distribution of $R_{i+1}$ is typically unknown, but the other component, i.e., the service time $B_{i+1}$, {\it is} known. 
This motivates our approximation in which $R_i$ is from the same class of distributions as $B_i$.

Now, starting with $r_1 = \beta_1$ and $v_1 = \sigma_1^2$, the core of our approximation algorithm is to iterate, for $i\in\{1,\ldots,n-1\}$, the following three steps:
\begin{description}
   \item[1.~Two-moment fit of $R_i$] Given the mean $r_i$ and variance $v_i$ corresponding to the sojourn time $R_{i}$, generate a two-moments fit of this sojourn time from one of the classes $ \{ {\rm Ph,~ W,~ LN} \}$.
   \item[2.~Compute ${\mathscr L}_{i+1} ( {x}_i)$] Determine the $(i+1)$-st summand of the objective function (i.e., ${\mathscr L}_{i+1} ( {x}_i)$)  
   using the two-moment fit of $R_i$. 
   \item[3.~Compute $r_{i+1}$ and $v_{i+1}$] Compute the mean $r_{i+1}$ and variance $v_{i+1}$ corresponding to the sojourn time $R_{i+1}$, by exploiting (\ref{Ri}) and (\ref{Vi}) and the two-moment fit of $R_i$. 
\end{description}
In Section~\ref{sec:alg} below, we will present closed-form expressions for the two-moments fits, $\ee(R_i - x)^+$, and the objects (\ref{Ri}) and (\ref{Vi}).
In addition, we detail the corresponding algorithms for the three different classes of two-moments fits. 

In our numerical experiments we evaluate the various approximations of the objective function $\objAS$ for given schedules ${\boldsymbol x}$. 
Furthermore, we assess the performance of corresponding {\it optimal} schedules by minimizing the approximations of $\objAS$ over ${\boldsymbol x}\in {\mathbb R}_+^{n-1}.$

\subsection{Comparison with the algorithm by De Kok \cite{de1989}}\label{COMP_DEK}
Our approach has borrowed elements from the approach proposed in the seminal paper \cite{de1989}, already briefly discussed in Section \ref{INTRO}. We substantially enrich the methodology developed in \cite{de1989}, though, in the following ways:
\begin{itemize}
    \item[$\circ$] Unlike \cite{de1989}, we have a focus on evaluating {\it transient} idle times and waiting times. More precisely, whereas \cite{de1989} aims to compute the {\it stationary} mean idle time and the {\it stationary} mean waiting time, we evaluate the {\it per-customer} mean idle time and mean waiting time. Due to \cite{de1989}'s inherent focus on stationarity it exclusively considers the situation of stochastically identical service times (i.e., the service times being independent and identically distributed), whereas our approach also covers heterogeneous customers (i.e., the service times being just independent). 
    \item[$\circ$] In \cite{de1989} the phase-type fit is suggested, without detailing the precise mapping of the parameters. Our work presents, in our setup that allows customer heterogeneity, this mapping via explicit expressions; see Algorithm \ref{algoPh} in Section~\ref{sec:alg}. 
    \item[$\circ$] We provide an in-depth analysis of alternative two-moments fits, namely those based on the Weibull and Lognormal distributions. We compare these to the one based on the phase-time fit that was relied upon in \cite{de1989}, considering both accuracy and execution times. 
    \item[$\circ$] We shed light on applying this type of approximation in the context of appointment scheduling. It means that we do not only consider the evaluation of the underlying objective function (in terms of accuracy and execution times), but also the evaluation of the optimizing appointment schedule. 
    In this application, we can streamline the procedure by working with the individual customer's {\it sojourn time} distributions. We perform a complexity analysis in terms of the number of customers $n$. 
\end{itemize} 
\section{Algorithms} \label{sec:alg}

In this section we present approximations for $\objAS$, $r_{i+1}$ and $v_{i+1}$, when $R_i$ is approximated by a low-dimensional phase-type random variable (Subsection~\ref{subs:PH}), a Weibull random variable (Subsection~\ref{subs:W}), and a Lognormal random variable (Subsection~\ref{subs:LN}).

\subsection{Phase-type fit} \label{subs:PH}
An idea that has been frequently applied in the literature, is that one can approximate non-negative random variables via phase-type distributions. 
Actually, this fit can be made arbitrarily accurate; see e.g.\ Asmussen \cite[Thm.\ III.4.2]{asmussen2003}. 
To obtain an accurate fit, however, the dimension of the underlying phase-type distribution can become large. 
This triggered the idea of phase-type approximations of relatively low dimension, chosen such that the approximating distribution has the same first two moments as the target distribution (the `two-moments fit'). 
In this low-dimensional approximation, the phase-type distributions used are mixtures of Erlang distributions and hyperexponential distributions. 
Due to the fact that phase-type distributions inherently lead to a Markovian system, explicit expressions for various queueing-related metrics are available, such as the idle and waiting time distributions (where in this context ‘explicit’ means in terms of  eigenvalues/eigenvectors of an associated eigensystem).

The two-moments fit we will be applying is the one presented in Tijms \cite{tijms1986}, matching the first and second moment, or, equivalently, the mean and the squared coefficient of variation ({\sc scv}). 
The {\sc scv} of the random variable $X$ is defined as ${\mathbb V}{\rm ar}\,X / ({\mathbb E}X)^2$. 
As prescribed by the fit presented in \cite{tijms1986}, we choose to match a mixture of two Erlang distributions in case the sojourn time $R_i$ has an {\sc scv} smaller than~1, and a hyperexponential distribution in case of an {\sc scv} larger than or equal to~1. 
We proceed by detailing this fitting procedure. 

\subsubsection*{Case $\mbox{\sc scv} \geqslant 1$}
We say that $R_i \sim  {\rm HE}(\alpha,\mu_1,\mu_2)$ in case the underlying density is defined through (for $y\geqslant 0$)
\[f_{R_i}(y) = \alpha \,\mu_1e^{-\mu_1 y} + (1-\alpha) \,\mu_2e^{-\mu_2 y} ,\]
for a probability $\alpha\in[0,1]$, and rates $\mu_1,\mu_2>0$.
It is readily verified that
\[{\mathbb E}\,R_i = \frac{\alpha}{\mu_1}+\frac{1-\alpha}{\mu_2},\:\:\:\:
{\mathbb V}{\rm ar}\,R_i =  \frac{\alpha}{\mu_1^2}+ \frac{1-\alpha}{\mu_2^2} + \alpha (1-\alpha) \Big(\frac{1}{\mu_1}-\frac{1}{\mu_2} \Big)^2.\]                                                     
For completeness, the two-moment fit is given in Appendix \ref{Explicit_FIT}.
To calculate the corresponding loss function, we obtain by an elementary calculation 
\begin{align} \label{eqn:excess-HE}
    {\mathbb E}(R_i-x)^+ & = \frac{\alpha}{\mu_1} e^{-\mu_1 x} + \frac{1-\alpha}{\mu_2}  e^{-\mu_2 x},
\end{align}
an expression that can also be derived relying on the memoryless property of the exponential distribution.
The above directly yields
\begin{equation}\label{DEFHE}
    {\mathscr L}_{i+1} ( {x})= \frac{\alpha}{\mu_1}(e^{-\mu_1 x} -\omega) + \frac{(1-\alpha)}{\mu_2}(e^{-\mu_2 x} -\omega) + \omega x =:{\mathscr L}^{\rm (HE)}(\alpha,\mu_1,\mu_2,x).
\end{equation}
Furthermore, using similar calculations combined with \eqref{Ri} and \eqref{Vi}, we find that
\begin{align} 
    r_{i+1}
    &= \frac{\alpha}{\mu_1} e^{-\mu_1 x_i}+\frac{1-\alpha}{\mu_2} e^{-\mu_2 x_i} +\beta_{i+1}  \label{RiHE} \\
    v_{i+1} & = 2\left(\frac{\alpha}{\mu_1^2} e^{-\mu_1 x_i}+\frac{1-\alpha}{\mu_2^2} e^{-\mu_2 x_i}\right) -\left( \frac{\alpha}{\mu_1} e^{-\mu_1 x_i}+\frac{1-\alpha}{\mu_2} e^{-\mu_2 x_i} \right)^2 +\sigma_{i+1}^2.   \label{ViHE}
\end{align}

\subsubsection*{Case $\mbox{\sc scv} < 1$}
Now suppose that $R_i$ is mixed Erlang, characterized by the probability $p\in[0,1]$, the shape parameter $k\in{\mathbb N}$ and the rate $\mu>0$, in the sequel denoted by $R_i \sim {\rm ME}(p,k,\mu)$. 
This concretely means that the underlying density is given by (for $y\geqslant 0$)
\begin{equation}\label{MED}
f_{R_i}(y) = p \,\frac{\mu^{k-1}y^{k-2}}{(k-2)!} e^{-\mu y}+ (1-p)\,\frac{\mu^ky^{k-1}}{(k-1)!} e^{-\mu y}.
\end{equation}
It requires a few elementary integrations to check that
\[{\mathbb E}\,R_i = p\frac{k-1}{\mu}+(1-p)\frac{k}{\mu}=\frac{k-1}{\mu}+\frac{1-p}{\mu},\:\:\:\:
{\mathbb V}{\rm ar}\,R_i = \frac{k-p^2}{\mu^2} .\]
While in Appendix~\ref{APP} we give alternative derivations for both $r_i$ and $v_i$ based on probabilistic arguments, in the expressions below we use direct (though lengthy) computations. 
More specifically, after a few steps of elementary calculus, involving a substitution, it follows that
\begin{align*}
\int_0^\infty y \mu^k (y+x)^{k-1} e^{-\mu (y+x)}\,{\rm d}y &=\frac{1}{\mu} \Gamma (k+1, \mu x) -x\, \Gamma(k, \mu x) \\
& = \Gamma(k, \mu x)\Big( \frac{k}{\mu}-x\Big)+ \mu^{k-1}x^k e^{-\mu x} =:M (k,x,\mu),
\end{align*}  
where $\Gamma(k,t)$ denotes the incomplete Gamma integral $\int_t^\infty z^{k-1} e^{-z}\,{\rm d}z$, which can be called in any standard numerical software package.
The above identity in combination with a few standard computations leads to  
\begin{align} 
    \notag {\mathbb E}(R_i-x)^+ 
    = \frac{p}{(k-2)!}\,M (k-1,x,\mu)+\frac{1-p}{(k-1)!}\,M (k,x,\mu). 
\end{align}
entailing that the corresponding summand in the loss function is
\begin{align} 
 {\mathscr L}_{i+1} ( {x}) &=  \frac{p}{(k-2)!}\,M (k-1,x,\mu)+\frac{1-p}{(k-1)!}\,M (k,x,\mu)  + \omega \left(x- \frac{k-p}{\mu}\right) 
    =:{\mathscr L}^{\rm (ME)}(p,k,\mu,x).
\label{DEFME}
\end{align}
Analogously, focusing on second moments, we find that
\begin{align*}
\int_0^\infty y^2 \mu^k &(y+x)^{k-1} e^{-\mu (y+x)}\,{\rm d}y  =\frac{1}{\mu^2}\Gamma(k+2, \mu x) - \frac{2x}{\mu} \Gamma(k+1, \mu x)+ x^2 \Gamma(k, \mu x) \\
& = \Gamma(k+1,\mu x)\Big( \frac{k+1-2\mu x}{\mu^2}\Big)+x^2\Gamma(k, \mu x)+ \mu^{k-1}x^{k+1}e^{-\mu x} \\
& = \Gamma(k,\mu x)\frac{k+(k-\mu x)^2}{\mu^2} + e^{-\mu x}\mu^{k-2}x^k(k+1-\mu x)  =:S (k,x,\mu).
\end{align*}
Hence, by applying \eqref{Ri} and \eqref{Vi},
\begin{align}
r_{i+1} &= \frac{p}{(k-2)!}\,M (k-1,x_i,\mu)+\frac{1-p}{(k-1)!}\,M (k,x_i,\mu)+\beta_{i+1} \label{RiME} \\
\notag v_{i+1}&=\frac{p}{(k-2)!}\,S (k-1,x_i,\mu)+\frac{1-p}{(k-1)!}\,S (k,x_i,\mu)\,\\ \label{ViME}
&\:\:\:-\left(\frac{p}{(k-2)!}\,M (k-1,x_i,\mu)+\frac{1-p}{(k-1)!}\,M (k,x_i,\mu)\right)^2+\sigma^2_{i+1}.
\end{align}

 


\vspace{3mm}

Using the concepts developed in Section \ref{sec:loss}, we are now in a position to describe our algorithm to approximate the objective function \eqref{OBJ}. 
The main idea is that for every client $i$ we approximate the mean $r_i$ and variance $v_i$ of her sojourn time, based on the mean $r_{i-1}$ and variance $v_{i-1}$ of the sojourn time of her predecessor, using the fit introduced above. 
The pseudocode for the resulting approximation algorithm, that we display below, is self-explanatory. 
Observe that the complexity of this algorithm (as well as those presented in the next subsections) is $O(n)$.

\vspace{3mm}

{\small
\begin{algorithm2e}[H]
\caption{Approximation algorithm: Phase-type fit}
 \label{algoPh}
 \DontPrintSemicolon
 \KwData{$(\beta_1,\ldots,\beta_n)$,  $(\sigma_1^2,\ldots,\sigma_n^2)$, $(x_1,\ldots,x_n)$.}
 \KwResult{$(r_1,\ldots,r_n)$,  $(v_1,\ldots,v_n),$ ${\mathscr L}^{({\rm Ph})} ( {\boldsymbol x})$.}
 Initialization: $i:=1$\;
 $r_1:=\beta_1$, $v_1:=\sigma_1^2$\;
  ${\mathscr L}^{({\rm Ph})} ( {\boldsymbol x}):=0$\;
 \While{$i < n$}{
  \eIf{$r_i^2\leqslant v_i$}{
   Perform HE phase-type fit: given $r_i$ and $v_i$, determine $\alpha,\mu_1,\mu_2$ using (\ref{eqn:H2-fit})\;
   ${\mathscr L}^{({\rm Ph})} ( {\boldsymbol x}):={\mathscr L}^{({\rm Ph})} ( {\boldsymbol x})+{\mathscr L}^{\rm (HE)}(\alpha,\mu_1,\mu_2,x_i)$, using the definition in \eqref{DEFHE}\;
    Determine $r_{i+1}$ by \eqref{RiHE}, $v_{i+1}$ by \eqref{ViHE}\;
   }{
   Perform ME phase-type fit: given $r_i$ and $v_i$, determine $p,k,\mu$ using (\ref{eqn:ME-fit})\;
   ${\mathscr L}^{({\rm Ph})} ( {\boldsymbol x}):={\mathscr L}^{({\rm Ph})} ( {\boldsymbol x})+{\mathscr L}^{\rm (ME)}(p,k,\mu,x_i)$, using the definition in \eqref{DEFME}\;
   Determine $r_{i+1}$ by \eqref{RiME}, $v_{i+1}$ by \eqref{ViME}\;
  }
  $i:=i+1$\;
 }
 
\end{algorithm2e}}

\vspace{3mm}

\subsection{Weibull fit} \label{subs:W}
We write $R_i \sim {\rm W}(\lambda,\alpha)$,  to denote that the density function underlying the sojourn time $R_i$ is given by (for $y\geqslant 0$)
\[
f_{R_i}(y)= \alpha \lambda^\alpha y^{\alpha-1} e^{-(\lambda y)^\alpha},
\]
for parameters $\lambda>0$ and $\alpha>0$; observe that this Weibull class includes distributions with tails that are both heavier and lighter tails than the tail of the exponential distribution.
It is well-known that
\[{\mathbb E}\,R_i= \lambda^{-1} \Gamma(1+1/\alpha),\:\:\:{\mathbb V}{\rm ar}\,R_i = \lambda^{-2}\big(\Gamma(1+2/\alpha)-(\Gamma(1+1/\alpha))^2\big).\]
As a consequence, the corresponding squared coefficient of variation equals
\[\frac{\Gamma(1+2/\alpha)-(\Gamma(1+1/\alpha))^2}{(\Gamma(1+1/\alpha))^2},\]
which decreases from $\infty$ to $0$ as $\alpha$ increases from $0$ to $\infty$, so that for each value of the squared coefficient of variation there is a unique $\alpha$.
This means that for each mean-variance pair there is a unique $\alpha$ and $\lambda$, which can be determined sequentially (first the $\alpha$, and then the $\lambda$, that is).

In order to evaluate the objective function we again first determine ${\mathbb E}(R_i-x)^+$. 
Applying integration by parts in the second equality, the substitution $w=(\lambda y)^\alpha$ in the fourth, and using the recurrence relation of the incomplete Gamma function, we obtain
\begin{align}
  {\mathbb E}(R_i-x)^+ \notag & = \int_x^\infty f_{R_i}(y) (y-x)\,{\rm d}y = \int_x^\infty (1-F_{R_i}(y))\,{\rm d}y= \int_x^\infty e^{-(\lambda y)^\alpha}\,{\rm d}y\\\label{EW}&=    \frac{1}{\lambda\alpha}\int_{(\lambda x)^\alpha}^\infty w^{1/\alpha -1} e^{-w}\,{\rm d}w = {\frac{1}{\lambda\alpha} \Gamma(1/\alpha, (\lambda x)^\alpha) } \\
  & = \frac{1}{\lambda} \Gamma(1+1/\alpha, (\lambda x)^\alpha) - x\,e^{-(\lambda x)^\alpha}.
\end{align}  
We thus find
\begin{align} 
    {\mathscr L}_{i+1} ( {x}) &=  
 \frac{1}{\lambda} \Gamma(1+1/\alpha, (\lambda x)^\alpha) - x\,e^{-(\lambda x)^\alpha} +\omega x -\frac{\omega}{\lambda}\Gamma(1+1/\alpha)=:{\mathscr L}^{\rm (W)}(\alpha,\lambda,x).
\label{DEW}
\end{align}
After some routine calculations, we obtain
\begin{align} 
    r_{i+1} &= \frac{1}{\lambda} \Gamma (1+{1}/{\alpha}, (\lambda x_i)^\alpha) - x_i e^{-(\lambda x_i)^\alpha} + \beta_{i+1}  \label{RiW} \\
\notag   
    v_{i+1} &= \frac{1}{\lambda^2}\Gamma(1+{2}/{\alpha},(\lambda x_i)^\alpha) -\frac{2x_i}{\lambda}\Gamma(1+{1}/{\alpha}, (\lambda x_i)^\alpha) + x_i^2 e^{-(\lambda x_i)^\alpha} \:\\ &\hspace{2cm}-\left(\frac{1}{\lambda} \Gamma (1+{1}/{\alpha}, (\lambda x_i)^\alpha) - x_i e^{-(\lambda x_i)^\alpha} \right)^2 + \sigma_{i+1}^2. \label{ViW}
\end{align}
Similarly to what we have done for the phase-type fit in Subsection~\ref{subs:PH}, the above can be summarized via pseudocode; see the algorithm below.

\vspace{3mm}

{\small
\begin{algorithm2e}[H]
\caption{Approximation algorithm: Weibull fit}
 \label{algoW}
\DontPrintSemicolon
 \KwData{$(\beta_1,\ldots,\beta_n)$,  $(\sigma_1^2,\ldots,\sigma_n^2)$, $(x_1,\ldots,x_n)$.}
 \KwResult{$(r_1,\ldots,r_n)$,  $(v_1,\ldots,v_n),$ ${\mathscr L}^{({\rm W})} ( {\boldsymbol x})$.}
 Initialization: $i:=1$\;
 $r_1:=\beta_1$, $v_1:=\sigma_1^2$\;
  ${\mathscr L}^{({\rm W})} ( {\boldsymbol x}):=0$\;
 \While{$i < n$}{
   Perform Weibull fit: given $r_i$ and $v_i$, determine $\lambda,\alpha$\;
   ${\mathscr L}^{({\rm W})} ( {\boldsymbol x}):={\mathscr L}^{({\rm W})} ( {\boldsymbol x})+{\mathscr L}^{\rm (W)}(\alpha,\lambda,x_i)$, using the definition in \eqref{DEW}\;
    Determine $r_{i+1}$ by \eqref{RiW}, $v_{i+1}$ by \eqref{ViW}\;
  $i:=i+1$\;
 }
 
\end{algorithm2e}}

\vspace{3mm}

\subsection{Lognormal fit} \label{subs:LN}
We write that $R_i$ is Lognormal, with parameters $\mu \in (-\infty, +\infty)$ and $\tau>0$, denoted by $R_i \sim {\rm LN}(\mu,\tau^2)$, if the underlying density is (for $y\geqslant 0$)
\[ { f_{R_i}(y)=\frac{1}{ \sqrt{2\pi}\tau y} \exp\left(-\frac{(\log y-\mu)^2}{2\tau^2}\right).}\]
In this case we have
\[{\mathbb E}\,R_i = e^{\mu+\tau^2/2},\:\:\:{\mathbb V}{\rm ar} \,R_i =(e^{\tau^2}-1) \,e^{2\mu+\tau^2}. \]
Observe that the corresponding squared coefficient of variation is $e^{\tau^2}-1$, which increases from~$0$ to $\infty$ as $\tau^2$ increases from~$0$ to $\infty$. 
Hence, to obtain a two-moments fit we first determine $\tau^2$ by equating the squared coefficient of variation to its target value, and then (given the value of $\tau^2$) the parameter $\mu$.

Applying the transformation $w=\log y$ in the second equality, and isolating the square in the third, 
\begin{align}\notag
  {\mathbb E}(R_i-x)^+ &  
  =  \int_x^\infty
  (y-x) \,\frac{1}{\sqrt{2\pi}\tau y} \exp\left(-\frac{(\log y-\mu)^2}{2\tau^2}\right)\,{\rm d}y\\\notag
  &=\int_{\log x}^\infty (e^w-x) \frac{1}{\sqrt{2\pi}\tau} \exp\left(-\frac{(w-\mu)^2}{2\tau^2}\right){\rm d}w \\\notag
  & = \frac{1}{\sqrt{2\pi}\tau}\left(e^{\mu+\tau^2/2}\int_{\log x}^\infty \exp\left(-\frac{(w-(\mu+\tau^2))^2}{2\tau^2}\right){\rm d}w - x\int_{\log x}^\infty \exp\left(-\frac{(w-\mu)^2}{2\tau^2}\right){\rm d}w 
  \right)\\
  &\notag = e^{\mu+\tau^2/2} \,{\mathbb P}({\rm N}(\mu+\tau^2,\tau^2)>\log x) - x \,{\mathbb P}({\rm N}(\mu,\tau^2)>\log x)\\
   &=e^{\mu+\tau^2/2} \Psi\left(\frac{\log x - \mu-\tau^2}{\tau}\right) - x \,\Psi\left(\frac{\log x-\mu}{\tau}\right),\label{ELN}
\end{align}
with $\Psi(x):=\int_x^\infty (\sqrt{2\pi})^{-1}\exp(-y^2/2)\,{\rm d}y$ denoting the standard Normal complementary distribution function, and $X\sim {\rm N}(\mu,\tau^2)$ denoting that $X$ is Normally distributed with mean $\mu$ and variance $\tau^2$. 
Hence, the objective function then reads
\begin{align} 
   {\mathscr L}_{i+1} ( {x}) &=  
   e^{\mu+\tau^2/2} \,\left(\Psi\left(\frac{\log x-\mu-\tau^2}{\tau}\right)-\omega\right) -  x\left(\Psi\left(\frac{\log x-\mu}{\tau}\right)-\omega\right)=:{\mathscr L}^{\rm (LN)}(\mu,\tau^2,x).
\label{DEFLN}
 \end{align}
Using similar steps as above to compute the variance, we have
\begin{align}
\label{RiLN} r_{i+1} &= e^{\mu+\tau^2/2} \,\Psi\left(\frac{\log x_i - \mu-\tau^2}{\tau}\right) - x_i \,\Psi\left(\frac{\log x_i -\mu }{\tau}\right) + \beta_{i+1} \\
\notag
v_{i+1} &=  e^{2(\mu + \tau^2)}\Psi\left( \frac{\log x_i-\mu -2\tau^2}{\tau} \right) -2x_i\, e^{\mu+\tau^2/2} \Psi\left(\frac{\log x_i - \mu-\tau^2}{\tau}\right) +x_i^2\, \Psi \left(\frac{\log x_i-\mu }{\tau} \right) \:\\
& \hspace{2cm} -\left(e^{\mu+\tau^2/2} \,\Psi\left(\frac{\log x_i - \mu-\tau^2}{\tau}\right) - x_i \,\Psi\left(\frac{\log x_i-\mu}{\tau}\right)\right)^2 +\sigma_{i+1}^2 .\label{ViLN}
\end{align}
The approach for the Lognormal case is summarized by the pseudocode below.

\vspace{3mm}

{\small 
\begin{algorithm2e}[H]
\caption{Approximation algorithm: Lognormal fit}
 \label{algoLN}
 \DontPrintSemicolon
 \KwData{$(\beta_1,\ldots,\beta_n)$,  $(\sigma_1^2,\ldots,\sigma_n^2)$, $(x_1,\ldots,x_n)$.}
 \KwResult{$(r_1,\ldots,r_n)$,  $(v_1,\ldots,v_n),$ ${\mathscr L}^{({\rm LN})} ( {\boldsymbol x})$.}
 Initialization: $i:=1$\;
 $r_1:=\beta_1$, $v_1:=\sigma_1^2$\;
  ${\mathscr L}^{({\rm LN})} ( {\boldsymbol x}):=0$\;
 \While{$i < n$}{
   Perform Lognormal fit: given $r_i$ and $v_i$, determine $\mu,\tau^2$ using (\ref{eqn:LN-fit})
   ${\mathscr L}^{({\rm LN})} ( {\boldsymbol x}):={\mathscr L}^{({\rm LN})} ( {\boldsymbol x})+{\mathscr L}^{\rm (LN)}(\mu,\tau^2,x_i)$, using the definition in \eqref{DEFLN}\;
    Determine $r_{i+1}$ by \eqref{RiLN}, $v_{i+1}$ by \eqref{ViLN}\;
  $i:=i+1$\;
 }
 
\end{algorithm2e}}

\vspace{3mm} 
\section{Numerical experimentation}\label{NUM}
This section presents an extensive numerical evaluation of the three approximation algorithms, in terms of accuracy (Subsections~\ref{EQD}--\ref{HET}) and efficiency (Subsection~\ref{comp_speed}). 
Most of the experiments consider situations with homogeneous customers, whereas we assess the impact of customer heterogeneity in one specific experiment in Subsection~\ref{HET}. For ease, we normalize time such that the mean service time equals 1.
The experiments corroborate the claim that our approximation method is remarkably accurate and of low computational complexity. 

We start by introducing convenient notation. 
The objective in appointment scheduling is to determine the optimal schedule, i.e., we wish to find the vector ${\boldsymbol x}^\star$ that minimizes the loss function:
\begin{equation}\label{minConv}{\boldsymbol x}^\star =\arg\min_{{\boldsymbol x}}{\mathscr L}({\boldsymbol x}).\end{equation}
As pointed out earlier in this paper, the objective function ${\mathscr L}({\boldsymbol x})$ can only be analytically evaluated for phase-type service times. 
In that case we have access to ${\mathscr L}({\boldsymbol x})$, so that the optimal schedule  ${\boldsymbol x}^\star$ can be determined. 
At the conceptual level the minimization in \eqref{minConv} is relatively easy, exploiting that ${\mathscr L}({\boldsymbol x})$ is convex in ${\boldsymbol x}$ \cite{kuiper2023}. 
A serious complication, however, lies in the numerical complexity: for large $n$ the evaluation of ${\boldsymbol x}^\star$ for phase-type service times is computationally demanding, as a consequence of the high dimension of matrices involved.
Recall from the introduction that, with $d_j$ denoting the dimension of the phase-type approximation of the service time of $j$-th client, the method requires various manipulations with matrices of dimension $D\times D$, with $D:=d_1+\cdots+d_n$, rendering the complexity of the evaluation of ${\mathscr L} ( {\boldsymbol x})$ to be $O(D^3)$. 
As we aim to find ${\boldsymbol x}^\star$, the function ${\mathscr L}({\boldsymbol x})$ has to be computed for a substantial number of schedules ${\boldsymbol x}$, making this approach infeasible. 
For this reason we wish to assess the performance of the three fast approximations that were proposed in the previous section.

As we need benchmarks to assess the accuracy of these approximations, we note that {\it for a given schedule} ${\boldsymbol x}$, one can still accurately approximate ${\mathscr L}({\boldsymbol x})$ using Monte Carlo techniques, in particular for non phase-type service times. 
However, finding the {\it optimal} schedule ${\boldsymbol x}^\star$ by estimating the objective function using Monte Carlo in each step of the minimization routine, is evidently not an efficient approach. 

\vb

Due to the intrinsic complications described above, various appointment scheduling heuristics have been proposed.
In the case of homogeneous clients, the most straightforward candidates are the {\it equidistant schedule} and the {\it Bailey-Welch rule}. 
In the equidistant rule the appointment schedule corresponds to constant interarrival times. In the sequel, we let 
${\boldsymbol x}_{\rm Eq}(y)$ denote such an equidistant schedule with interarrival times equal to $y>0$; given the mean service time is 1, one typically chooses values of $y$ larger than 1.
In the Bailey-Welch rule, the two first customers are scheduled to arrive at time $0$, and the $i$-th (for $i=3,4,\ldots$) at time $i-2$ (recalling that we put the mean service time equal to 1); we denote the corresponding schedule by ${\boldsymbol x}_{\rm BW}.$ 
It is anticipated that for somewhat larger numbers of customers, the standard Bailey-Welch rule may perform poorly: apart from the way the first customer is treated, it effectively corresponds to a critically loaded queue (i.e., a queue with load $1$) in which, in particular at the end of the schedule, waiting times may become substantial. 
To remedy this undesired effect, one could think of a hybrid form of the equidistant schedule and the Bailey-Welch rule ${\boldsymbol x}_{\rm BW}(y)$: the two first customers are scheduled to arrive at time $0$, and the $i$-th (for $i=3,4,\ldots$) at time $(i-2)y$ for some $y$ larger than 1. 

In the remainder of this section we discuss the results of our numerical experiments. 
In all experiments that we present in this section, unless stated otherwise, we chose 
\[\omega = 0.5,\:\:\:n=40,\:\:\:\:\beta_1=\cdots=\beta_n=1,\:\:\:\:\sigma^2:=\sigma_1^2=\cdots=\sigma_n^2 \in\{0.4,0.7,1.0,1.3\},\] 
but we remark that additional experiments reveal that the below observations carry over to instances with $\omega \in \{0.25, 0.75\}$ and $n=10$. 
If the benchmark value of the objective function has been obtained via Monte Carlo simulation, then it was estimated based on $10^5$ runs (ensuring a sufficiently high precision); we then use the notation ${\mathscr L}_{\rm sim}$.
All numerical experiments have been done using an {\sc matlab} R2022a implementation, on a standard MacBook Pro.

\subsection{Equidistant schedules}\label{EQD}
We start by analyzing the accuracy of our approximations for equidistant schedules, considering the cases $y\in\{1.2,1.5,1.8\}$. 
We assess the performance of the approximations based on relative errors.
\begin{itemize}
    \item[$\circ$]
For phase-type service times, we wish to compare the `slow' objective function ${\mathscr L}({\boldsymbol x}_{\rm Eq}(y))$ with its `fast' counterpart ${\mathscr L}^{(\rm Ph)}({\boldsymbol x}_{\rm Eq}(y))$. 
We introduce, for ease leaving out the argument ${\boldsymbol x}_{\rm Eq}(y)$,
\[ \Delta_{\rm Ph}^{\rm ({Ph})}= \Bigg| \frac{{\mathscr L}^{\rm ({Ph})}-{\mathscr L}}{{\mathscr L}} \Bigg| \cdot 100\% .\]
 \item[$\circ$] 
For the case of Weibull service times, we wish to compare ${\mathscr L}({\boldsymbol x}_{\rm Eq}(y))$ (which is estimated through Monte Carlo simulation) with the fast approximations ${\mathscr L}^{(\rm W)}({\boldsymbol x}_{\rm Eq}(y))$ and ${\mathscr L}^{(\rm Ph)}({\boldsymbol x}_{\rm Eq}(y))$.
Likewise, for Lognormal service times we compare ${\mathscr L}({\boldsymbol x}_{\rm Eq}(y))$ (again estimated through Monte Carlo simulation) with ${\mathscr L}^{(\rm LN)}({\boldsymbol x}_{\rm Eq}(y))$ and ${\mathscr L}^{(\rm Ph)}({\boldsymbol x}_{\rm Eq}(y))$.
We now work with, again omitting the argument ${\boldsymbol x}_{\rm Eq}(y)$,
\[ \Delta_{\rm C}^{\rm (R)}=\Bigg| \frac{{\mathscr L}^{\rm (R)}-{\mathscr L}_{\rm sim}}{{\mathscr L}_{\rm sim}}\Bigg|\cdot 100\%,\]
where R represents the type of distribution that has been used to fit the sojourn times, and 
${\rm C}\in \{ {\rm W,~ LN}  \}$ is the class the service times stem from (under which the simulation has been performed).
\end{itemize}



Table \ref{equi} presents the actual and approximating objective functions, as well as the relative errors, for phase-type, Weibull, and Lognormal service time distributions.
The main conclusion from this table is that, across all instances, the fast phase-type approximation ${\mathscr L}^{(\rm Ph)}({\boldsymbol x}_{\rm Eq}(y))$ accurately approximates the objective function ${\mathscr L}({\boldsymbol x}_{\rm Eq}(y))$. 
One could have anticipated that if the service times have a Weibull (Lognormal, respectively) distribution, then the Weibull approximation ${\mathscr L}^{(\rm W)}({\boldsymbol x}_{\rm Eq}(y))$ (Lognormal approximation ${\mathscr L}^{(\rm LN)}({\boldsymbol x}_{\rm Eq}(y))$, respectively) would have worked significantly better, but the table does not provide convincing support for such a claim; it is true that the error is typically slightly smaller, but the differences tend to be minimal.
For Lognormal service times, the approximation ${\mathscr L}^{(\rm LN)}({\boldsymbol x}_{\rm Eq}(y))$ performs actually somewhat worse than ${\mathscr L}^{(\rm Ph)}({\boldsymbol x}_{\rm Eq}(y))$ in the regime that the {\sc scv} is relatively large and $y$ is relatively low, an effect that may be explained by the heavy tail of the Lognormal distribution.

\begin{table}[h]
{\small
\begin{tabular}{ll lll |lll| lll| lll}
\toprule
\multicolumn{2}{l}{\textbf{\sc scv}} & \multicolumn{3}{c}{0.4} & \multicolumn{3}{c}{0.7} & \multicolumn{3}{c}{1} & \multicolumn{3}{c}{1.3} \\
\multicolumn{2}{l}{$\mathbf{\it y}$} & 1.2 & 1.5 & 1.8 & 1.2 & 1.5 & 1.8 & 1.2 & 1.5 & 1.8 & 1.2 & 1.5 & 1.8 \\
\midrule
\multirow{3}{*}{Phase-type} & ${\mathscr L}$ & 17.13 & 14.02 & 17.66 & 26.43& 18.55& 20.17   &34.53 &23.40 & 23.19  &41.26 & 28.14&26.45 \\
 & ${\mathscr L}^{(\rm Ph)}$  & 17.15 & 13.95 & 17.63 & 26.57 & 18.50 & 20.13  & 34.37 & 23.39 & 23.19   & 39.83 & 27.78 & 26.43\\
& $\Delta_{\rm Ph}^{(\rm Ph)}$ & 0.12 & 0.50 & 0.17& 0.56 &0.32 & 0.25& 0.49& 0.00 & 0.00& 3.56 & 1.33& 0.08\\
\midrule
\multirow{4}{*}{Weibull} & ${\mathscr L}_{\rm sim}$& 17.04 & 13.83 & 17.49 & 26.53 &18.59 & 20.18 & 34.61 & 23.41 & 23.19 & 41.63 & 28.10& 26.30\\
& ${\mathscr L}^{(\rm W)}$& 17.10& 13.76& 17.46& 26.59 & 18.56 & 20.16 & 34.24 & 23.41 & 23.19 & 40.59 & 28.04 & 26.31 \\
& $\Delta_{\rm W}^{\rm (W)}$& 0.35 &0.51 & 0.17& 0.19 & 0.16& 0.10 & 1.10& 0.04 & 0.00 & 2.50 & 0.21 & 0.08 \\
&$\Delta_{\rm W}^{\rm (Ph)}$ & 0.65 & 0.87 & 0.80 & 0.15 & 0.48 & 0.30 & 0.72 & 0.09 & 0.00 & 4.32 & 1.14 & 0.53\\
\midrule
\multirow{4}{*}{Lognormal} & ${\mathscr L}_{\rm sim}$ & 17.40& 14.53 & 18.11 & 26.22 & 19.23 & 20.93 & 33.49 & 23.80 & 23.91 & 39.62 & 28.08 & 26.88\\
& ${\mathscr L}^{(\rm LN)}$ & 17.00 & 14.46 & 18.07 & 24.55 & 19.08 & 20.87 & 30.13 & 23.33 & 23.79 & 34.54 & 27.03 & 26.62\\
& $\Delta^{\rm (LN)}_{\rm LN}$ & 2.50 & 0.48& 0.22& 6.37 & 0.83 & 0.33 & 10.06 & 2.02 & 0.50 & 12.85 & 3.78 & 0.97 \\
& $\Delta_{\rm LN}^{\rm (Ph)}$& 1.44 & 3.99 & 2.65 & 1.33 & 3.85 & 3.87 & 2.60 & 1.72 & 3.01& 0.53 & 1.10 & 1.67\\
\bottomrule
\end{tabular}
\caption{\label{equi} Performance assessment for equidistant schedules.  For all the loss functions in this table the schedule ${\boldsymbol x}_{\rm Eq}(y)$ has been used, with $y\in\{1.2,1.5,1.8\}$. }
}
\end{table}

\subsection{Bailey-Welch schedules\label{Bw_section}} 
In this case, we focus on $y=1.2$. Larger values of $y$ are less interesting, as then the equidistant schedules and those generated by the Bailey-Welch rule virtually coincide. 
We use a similar comparison and notation as in Subsection \ref{EQD}.
The conclusions from Table \ref{baileyW} are similar to those drawn from Table \ref{equi}. 
In particular, for all three types of service-time distributions, the fast phase-type approximation ${\mathscr L}^{(\rm Ph)}({\boldsymbol x}_{\rm BW})$  accurately approximates the objective function ${\mathscr L}({\boldsymbol x}_{\rm BW}(y))$. 
Again, for Lognormal service times the performance of ${\mathscr L}^{(\rm LN)}({\boldsymbol x}_{\rm BW}(y))$ is relatively poor, especially for larger values of the {\sc scv} (where it plays a role that the value of $y$ that we have chosen is relatively small).

\begin{table}[h]
\begin{tabular}{llllll}
\toprule 
\multicolumn{2}{l}{\textbf{\sc scv}} & 0.4 & 0.7 & 1 & 1.3 \\
\midrule
\multirow{3}{*}{Phase-type} & ${\mathscr L}$ & 18.74 & 28.11 & 36.17 & 42.81\\
& ${\mathscr L}^{(\rm Ph)}$ &18.78 & 28.27& 35.99 & 41.33\\
& $\Delta^{\rm (Ph)}_{\rm Ph}$ & 0.21 & 0.59 & 0.49 & 3.46\\
\midrule
\multirow{4}{*}{Weibull} & ${\mathscr L}_{\rm sim}$ & 18.70 & 28.21 & 36.25 & 43.21\\
& ${\mathscr L}^{(\rm W)}$ & 18.77 & 28.28 & 35.85 & 42.12\\
&$\Delta^{\rm (W)}_{\rm W}$ & 0.37 & 0.24 & 1.10 & 2.52\\
& $\Delta_{\rm W}^{\rm (Ph)}$ & 0.44 & 0.23 & 0.72 & 4.55\\
\midrule
\multirow{4}{*}{Lognormal} & ${\mathscr L}_{\rm sim}$ & 18.90& 27.74 & 34.97& 41.06\\
& ${\mathscr L}^{(\rm LN)}$ & 18.47 & 25.93 & 31.39 & 35.70\\
& $\Delta^{(\rm LN)}_{\rm LN}$ & 2.27 & 6.51 & 10.22 & 13.03 \\
& $\Delta_{\rm LN}^{\rm (Ph)}$ & 0.64 & 1.89 & 2.84 & 0.65\\
\bottomrule
\end{tabular}
\caption{\label{baileyW} Performance assessment for Bailey-Welch schedules. For all the loss functions in this table the schedule ${\boldsymbol x}_{\rm BW}(y)$ has been used, with $y=1.2$.}
\end{table}

\subsection{Optimized schedules}\label{optim_sc}
In this case, we study the accuracy of our approximations for optimized interarrival times. 
As indicated before, only for phase-type service times we can evaluate the objective function, and even then finding the optimal schedule $x^\star$ is computationally hard. The obvious alternative is to produce optimized interarrival times based on the approximate objective functions. To this end, 
the following notation is used:
\[{\boldsymbol x}_{\rm Ph}^\star =\arg\min_{{\boldsymbol x}}{\mathscr L}^{({\rm Ph})}({\boldsymbol x}),\:\:\:{\boldsymbol x}^\star_{\rm W} =\arg\min_{{\boldsymbol x}}{\mathscr L}^{({\rm W})}({\boldsymbol x}),\:\:\:{\boldsymbol x}^\star_{\rm LN} =\arg\min_{{\boldsymbol x}}{\mathscr L}^{({\rm LN})}({\boldsymbol x}).\]
\begin{itemize}
    \item[$\circ$]
For phase-type service times, our goal is to compare the value of our `fast' objective function ${\mathscr L}^{(\rm Ph)}$ evaluated at the corresponding optimized interarrival times ${\boldsymbol x}_{\rm Ph}^\star$, with the `slow' counterparts ${\mathscr L}({\boldsymbol x}^\star)$  and ${\mathscr L}({\boldsymbol x}_{\rm Ph}^\star)$. 
 \item[$\circ$] 
For the case of Weibull service times, we evaluate ${\mathscr L}^{({\rm W})}({\boldsymbol x}_{\rm W}^\star)$,  ${\mathscr L}({\boldsymbol x}_{\rm W}^\star)$, and ${\mathscr L}({\boldsymbol x}_{\rm Ph}^\star)$, where it is noted that ${\mathscr L}^{({\rm W})}({\boldsymbol x}_{\rm W}^\star)$ allows quick evaluation, while the evaluation of ${\mathscr L}({\boldsymbol x}_{\rm W}^\star)$, and ${\mathscr L}({\boldsymbol x}_{\rm Ph}^\star)$ is considerably more time-consuming.
As pointed out above, both ${\mathscr L}({\boldsymbol x}_{\rm W}^\star)$ and ${\mathscr L}({\boldsymbol x}_{\rm Ph}^\star)$ are estimated by Monte Carlo simulation.
Here, ${\boldsymbol x}_{\rm W}^\star$ and ${\boldsymbol x}_{\rm Ph}^\star$ are the optimal interarrival times based on approximate objective functions, where ${\boldsymbol x}_{\rm Ph}^\star$ is of specific interest in view of the computational speed (see Subsection~\ref{comp_speed} below).
 \item[$\circ$] 
Similar to above, for Lognormal service times, we compare  ${\mathscr L}^{({\rm LN})}({\boldsymbol x}_{\rm LN}^\star)$, ${\mathscr L}({\boldsymbol x}_{\rm LN}^\star)$, and ${\mathscr L}({\boldsymbol x}_{\rm Ph}^\star)$, where the latter two quantities are estimated by Monte Carlo simulation. 
\end{itemize}

Along with the relative error defined above, we also compute an `optimality gap'.
For phase-type service time, this gap is defined as 
\[ \Delta_{\rm opt}=\Bigg|\frac{{\mathscr L}({\boldsymbol x}_{\rm Ph}^\star)-{\mathscr L}({\boldsymbol x}^\star)}{{\mathscr L}({\boldsymbol x}^\star)}\Bigg| \cdot 100\% ,\]
whereas for Weibull and Lognormal service times we compute
\[ \Tilde{\Delta}_{\rm opt}=\Bigg|\frac{{\mathscr L}_{\rm sim}({\boldsymbol x}_{\rm C}^\star)-{\mathscr L}_{\rm sim}({\boldsymbol x}_{\rm Ph}^\star)}{\min\big\{ {\mathscr L}_{\rm sim}({\boldsymbol x}_{\rm C}^\star), ~{\mathscr L}_{\rm sim}({\boldsymbol x}_{\rm Ph}^\star)\big\} }\Bigg| \cdot 100\%,\]
with ${\rm C}\in\{{\rm W, LN}\}.$

Table \ref{optimSched} reveals that the Weibull and Lognormal approximations perform essentially equally well as the phase-type approximation (even for larger values of {\sc scv}, that is; for small {\sc scv} and lognormal service times, the phase-type approximation slightly deteriorates).
The consequence of this observation is that, when the objective is to produce optimal appointment schedules, we can use any of the three approximations. 
This makes the question relevant which of the three approximations requires the lowest computational effort; this question will be extensively studied in Section \ref{comp_speed}.

\begin{table}[h]
\begin{tabular}{llllll}
\toprule 
\multicolumn{2}{l}{\textbf{\sc scv}} & 0.4 & 0.7 & 1 & 1.3 \\
\midrule
\multirow{5}{*}{Phase-type} & ${\mathscr L}({\boldsymbol x}^\star )$   &13.59 & 18.37& 22.45 & 26.09 \\ 
& ${\mathscr L}({\boldsymbol x}^\star _{\rm Ph})$ &  13.61 &18.40 &22.53 & 26.21 \\ 
& $\Delta_{\rm opt}$ & 0.08 & 0.21 & 0.36 & 0.46 \\ \cline{2-6}
&  ${\mathscr L}^{(\rm Ph)}({\boldsymbol x}^\star _{\rm Ph})$  &13.52 &18.31 &22.45 &26.03 \\
& $\Delta^{\rm (Ph)}_{\rm Ph}$ & 0.65& 0.52& 0.35 & 0.69\\ 

\midrule

\multirow{6}{*}{Weibull} & ${\mathscr L}_{\rm sim}({\boldsymbol x}^\star _{\rm W})$ & 13.65& 18.44 &22.55 & 26.20 \\ 
& ${\mathscr L}_{\rm sim}({\boldsymbol x}^\star _{\rm Ph})$  & 13.47  & 18.47& 22.55 & 26.28\\ 
& $\Tilde{\Delta}_{\rm opt}$ & 1.35 & 0.12 & 0.00 & 0.30\\ \cline{2-6}
& ${\mathscr L}^{(\rm W)}({\boldsymbol x}^\star_{\rm Ph})$ & 13.35 & 18.36 & 22.45 & 25.98 \\
&$\Delta^{\rm (W)}_{\rm W}$ & 0.96 & 0.55 & 0.46 & 1.16 \\
&$\Delta^{\rm (Ph)}_{\rm W}$ & 0.33 & 0.85 & 0.46 & 0.97\\
\midrule

\multirow{6}{*}{Lognormal} & ${\mathscr L}_{\rm sim}({\boldsymbol x}^\star _{\rm LN})$ & 13.66 & 18.47 & 22.61 & 26.30\\ 
& ${\mathscr L}_{\rm sim}({\boldsymbol x}^\star _{\rm Ph})$ & 14.14 & 19.12 & 23.12 & 26.52\\ 
& $\Tilde{\Delta}_{\rm opt}$ & 3.57 & 3.49 & 2.25 & 0.87 \\ \cline{2-6}
&${\mathscr L}^{(\rm LN)}({\boldsymbol x}^\star _{\rm Ph})$ & 13.98 & 18.89 & 22.79 & 26.03 \\ 
&$\Delta^{\rm (LN)}_{\rm LN}$ & 1.15 & 1.19 & 1.45 & 1.87 \\
&$\Delta^{\rm (Ph)}_{\rm LN}$ & 4.43 & 4.23 & 2.92 & 1.88 \\

\bottomrule
\end{tabular}
\caption{\label{optimSched} Performance assessment for optimized schedules. }
\end{table}

\subsection{Customer heterogeneity}\label{HET}
To verify to what extent the above conclusions carry over to scenarios with heterogeneous service times, we performed a series of additional experiments. 
There are clearly many ways to introduce heterogeneity; the specific class of instances we pick in this experiment produces results that are representative for what we found in other experiments. 
Again, we take $\beta_1=\cdots=\beta_n=1$, but now some customers have $\sigma^2=0.7$ and some have $\sigma^2=1.3$. 
We consider the a set of scenarios corresponding to 40 customers that are arranged in eight batches of 5, where the {\sc scv} is kept fixed within each batch. 
The following list, with $P_1$ corresponding to an {\sc scv} of $0.7$, and $P_2$ to an {\sc scv} of $1.3$, encodes the instances:
\begin{itemize}
    \item[(A)]~$P_1$~---~$P_1$~---~$P_1$~---~$P_1$~---~$P_2$~---~$P_2$~---~$P_2$~---~$P_2$, 
    \item[(B)]~$P_2$~---~$P_2$~---~$P_2$~---~$P_2$~---~$P_1$~---~$P_1$~---~$P_1$~---~$P_1$,
    \item[(C)]~$P_1$~---~$P_1$~---~$P_2$~---~$P_2$~---~$P_1$~---~$P_1$~---~$P_2$~---~$P_2$,
    \item[(D)]~$P_2$~---~$P_2$~---~$P_1$~---~$P_1$~---~$P_2$~---~$P_2$~---~$P_1$~---~$P_1$,
    \item[(E)]~$P_1$~---~$P_2$~---~$P_1$~---~$P_2$~---~$P_1$~---~$P_2$~---~$P_1$~---~$P_2$,
    \item[(F)]~$P_2$~---~$P_1$~---~$P_2$~---~$P_1$~---~$P_2$~---~$P_1$~---~$P_2$~---~$P_1$.
\end{itemize}

In Table \ref{heterogeneous} we take an equidistant schedule with $y=1.5$, whereas we consider optimized interarrival times (as in Subsection~\ref{optim_sc}) again in Table \ref{het_opt_n}.
The results in Tables \ref{heterogeneous} and \ref{het_opt_n} confirm what was observed in our earlier experiments.
Most importantly, across the board the phase-type approximation ${\mathscr L}^{(\rm Ph)}({\boldsymbol x}_{\rm Eq})$ performs well (but better in the Weibull case than in the Lognormal case).

\begin{table}[h]
\begin{tabular}{llllllll}
\toprule
\multicolumn{2}{l}{\textrm{Scenario}} & A & B & C & D & E & F \\
\midrule
\multirow{3}{*}{Phase-type} &  ${\mathscr L}$ & 22.69 & 24.05 & 22.93 & 23.86 & 23.14 & 23.66 \\
   & ${\mathscr L}^{(\rm Ph)}$ & 22.40 & 24.07 & 22.78 & 23.89 & 23.07 & 23.68 \\
   & $\Delta^{\rm (Ph)}_{\rm Ph}$ & 1.30 & 0.06 & 0.66 & 0.16 & 0.30 & 0.08 \\ \midrule
\multirow{4}{*}{Weibull} & ${\mathscr L}_{\rm sim}$   & 22.70 & 23.99 & 22.90 & 23.78 & 23.08 & 23.59 \\
   & ${\mathscr L}^{(\rm W)}$ & 22.64 & 24.00 & 22.88 & 23.86 & 23.12 & 23.67  \\
   &$\Delta^{\rm (W)}_{\rm W}$ & 0.26 & 0.04 & 0.07 & 0.33 & 0.16 & 0.33  \\
   &$\Delta^{\rm (Ph)}_{\rm W}$ & 1.32 & 0.33 & 0.53 & 0.48 & 0.02 & 0.40  \\\midrule
\multirow{4}{*}{Lognormal} & ${\mathscr L}_{\rm sim}$  & 22.88 & 24.48 & 23.22 & 24.22 & 23.45 & 23.99  \\
   & ${\mathscr L}^{(\rm LN)}$ & 22.17 & 24.30 & 22.70 & 23.90 & 23.02 & 23.63 \\
   & $\Delta^{\rm (LN)}_{\rm LN}$& 3.10 & 0.74 & 2.22 & 1.30 & 1.85 & 1.51 \\
   &$\Delta^{\rm (Ph)}_{\rm LN}$ & 2.10 & 1.67 & 1.90 & 1.34 & 1.60 & 1.28 \\ \bottomrule
\end{tabular}
\caption{\label{heterogeneous} Performance assessment for heterogeneous scenarios.  For all the loss functions in this table the schedule ${\boldsymbol x}_{\rm Eq}(y)$ has been used, with $y=1.5$.}
\end{table}

\begin{table}[h]
\begin{tabular}{llllllll}
\toprule
\multicolumn{2}{l}{\textrm{Scenario}} & A & B & C & D & E & F \\
\midrule
\multirow{4}{*}{Phase-type} & ${\mathscr L}({\boldsymbol x}^\star )$ & 21.97 & 22.56 & 22.07 & 22.60 & 22.20 & 22.61  \\
& ${\mathscr L}({\boldsymbol x}^\star_{\rm Ph} )$ & 21.98 & 22.57& 22.08 &22.62 & 22.21& 22.64\\
   & ${\mathscr L}^{\rm (Ph)}({\boldsymbol x}^\star _{\rm Ph})$ & 21.84 & 22.55  & 21.94  & 22.59  & 22.12  & 22.61 \\
   & $\Delta_{\rm opt}$ & 0.04 & 0.05 &0.07 &0.10 & 0.04& 0.11\\
   & $\Delta^{\rm (Ph)}_{\rm Ph}$ & 0.60 & 0.06 & 0.58 & 0.05 & 0.38 & 0.00\\ \midrule

\multirow{4}{*}{Weibull} & ${\mathscr L}_{\rm sim}({\boldsymbol x}^\star _{\rm W})$ & 22.02 & 22.62 & 22.05 & 22.55 & 22.16 & 22.62 \\
   & ${\mathscr L}_{\rm sim}({\boldsymbol x}^\star _{\rm Ph})$ & 22.03 & 22.58 & 22.35 & 22.61 & 22.21 & 22.62 \\
   & $\Tilde{\Delta}_{\rm opt}$ & 0.06 & 0.17 & 1.36 & 0.27 & 0.23 & 0.09\\
   &$\Delta^{\rm (Ph)}_{\rm W}$ & 0.87 & 0.14 & 1.83 & 0.09 & 0.41 & 0.04 \\ \midrule
\multirow{4}{*}{Lognormal} & ${\mathscr L}_{\rm sim}({\boldsymbol x}^\star _{\rm LN})$ & 22.46 & 23.30 & 22.60 & 23.27 & 22.76 & 23.23 \\
   & ${\mathscr L}_{\rm sim}({\boldsymbol x}^\star _{\rm Ph})$ & 22.41 & 23.33 & 22.58 & 23.31 & 22.76 & 23.24\\
   & $\Tilde{\Delta}_{\rm opt}$ & 0.22 & 0.13 & 0.09 &0.17 & 0.00 & 0.04\\
   &$\Delta^{\rm (Ph)}_{\rm LN}$ & 2.54 & 3.34 & 2.83  & 3.09 & 2.81  & 2.70\\ \bottomrule
\end{tabular}
\caption{\label{het_opt_n} Performance assessment for heterogeneous scenarios with optimized interarrival times.}
\end{table}
\subsection{Computation speed}\label{comp_speed}
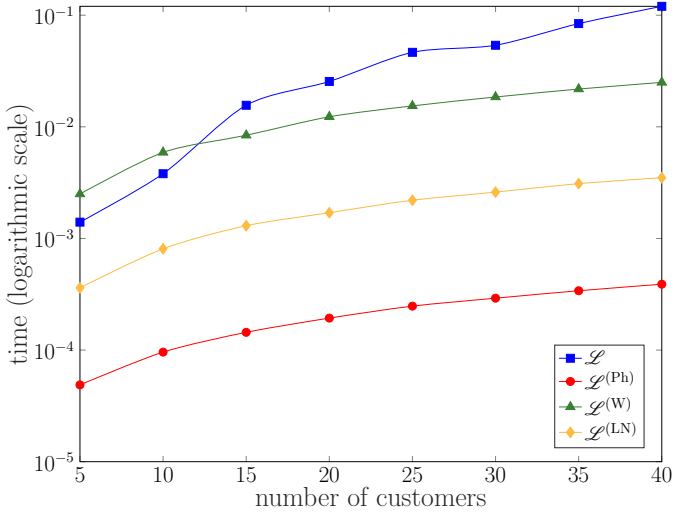
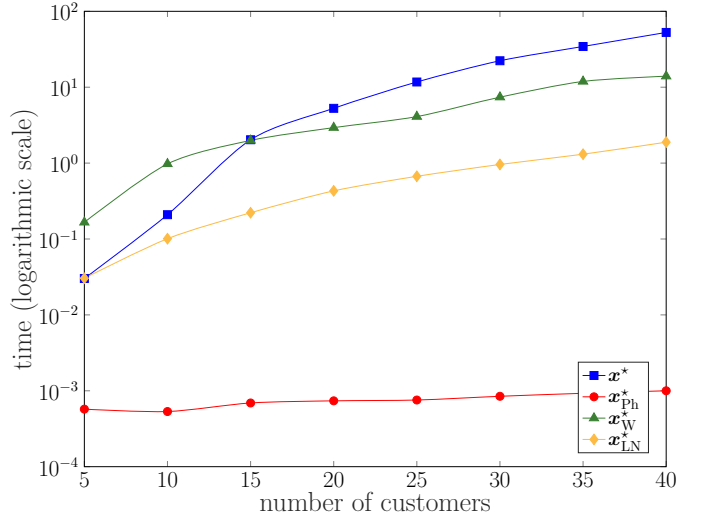
\begin{figure}[htbp]
    \begin{subfigure}{0.5\textwidth} 
        \centering
       \definecolor{mycolor1}{rgb}{1.00000,0.00000,1.00000}%
\begin{tikzpicture}[scale=0.5]

\begin{axis}[%
width=6.028in,
height=4.754in,
at={(1.011in,0.642in)},
scale only axis,
xmin=5,
xmax=40,
xtick={0,5,10,15,20,25,30,35,40},
xticklabel style={font=\fontsize{16}{14}\selectfont},
xlabel style={font=\color{white!15!black}},
xlabel={\fontsize{22}{14}\selectfont number of customers},
ymode=log,
ymin=1e-05,
ymax=0.1201,
yminorticks=false,
yticklabel style={font=\fontsize{16}{14}\selectfont},
ylabel style={font=\color{white!15!black}, yshift=10pt},
ylabel={\fontsize{22}{14}\selectfont time (logarithmic scale)},
axis background/.style={fill=white},
legend style={at={(0.97,0.03)}, anchor=south east, legend cell align=left, align=left, draw=white!15!black, font=\fontsize{15}{14}\selectfont}
]
\addplot [smooth, mark=square*, blue, mark size=3pt]
  table[row sep=crcr]{%
5	0.0014\\
10	0.0038\\
15	0.0156\\
20	0.0255\\
25	0.0465\\
30	0.0538\\
35	0.0841\\
40	0.1201\\
};
\addlegendentry{$\mathscr{L}$}

\addplot [smooth,mark=*,red, mark size=3pt]
  table[row sep=crcr]{%
5	4.8757e-05\\
10	9.586e-05\\
15	0.00014417\\
20	0.00019334\\
25	0.00024748\\
30	0.0002915\\
35	0.00034007\\
40	0.00038893\\
};
\addlegendentry{$\mathscr{L}^{(\rm Ph)}$}

\addplot [smooth, mark=triangle*, OliveGreen, mark size=4pt]
  table[row sep=crcr]{%
5	0.0025\\
10	0.0059\\
15	0.0084\\
20	0.0123\\
25	0.0154\\
30	0.0185\\
35	0.0218\\
40	0.025\\
};
\addlegendentry{$\mathscr{L}^{(\rm W)}$}

\addplot [smooth, mark=diamond*, Dandelion, mark size=4pt]
  table[row sep=crcr]{%
5	0.00036105\\
10	0.0008076\\
15	0.0013\\
20	0.0017\\
25	0.0022\\
30	0.0026\\
35	0.0031\\
40	0.0035\\
};
\addlegendentry{$\mathscr{L}^{(\rm LN)}$}

\end{axis}

\end{tikzpicture}%
        \caption{Loss functions}
        \label{fig:subfig1a}
    \end{subfigure}%
    \begin{subfigure}{0.5\textwidth} 
        \centering
       \definecolor{mycolor1}{rgb}{1.00000,0.00000,1.00000}%
\begin{tikzpicture}[scale=0.5]

\begin{axis}[%
width=6.028in,
height=4.754in,
at={(1.011in,0.642in)},
scale only axis,
xmin=5,
xmax=40,
xtick={0,5,10,15,20,25,30,35,40},
xticklabel style={font=\fontsize{16}{14}\selectfont},
xlabel style={font=\color{white!15!black}},
xlabel={\fontsize{22}{14}\selectfont number of customers},
ymode=log,
ymin=1e-04,
ymax=100,
yminorticks=false,
yticklabel style={font=\fontsize{16}{14}\selectfont},
ylabel style={font=\color{white!15!black}, yshift=10pt},
ylabel={\fontsize{22}{14}\selectfont time (logarithmic scale)},
axis background/.style={fill=white},
legend style={at={(0.97,0.03)}, anchor=south east, legend cell align=left, align=left, draw=white!15!black, font=\fontsize{15}{14}\selectfont}
]
\addplot [smooth, mark=square*, blue, mark size=3pt]
  table[row sep=crcr]{%
5	0.0302\\
10	0.209\\
15	2.0307\\
20	5.2713\\
25	11.7156\\
30	22.3102\\
35	34.4167\\
40	52.6915\\
};
\addlegendentry{${\boldsymbol x}^\star$}

\addplot [smooth,mark=*,red, mark size=3pt]
  table[row sep=crcr]{%
5	0.00057311\\
10	0.00053298\\
15	0.00069113\\
20	0.00073734\\
25	0.00075564\\
30	0.00084637\\
35	0.00092662\\
40	0.001\\
};
\addlegendentry{${ \boldsymbol x}^\star_{\rm Ph}$}

\addplot [smooth, mark=triangle*, OliveGreen, mark size=4pt]
  table[row sep=crcr]{%
5	0.1658\\
10	0.9765\\
15	1.9912\\
20	2.9342\\
25	4.1027\\
30	7.3928\\
35	11.9174\\
40	13.9645\\
};
\addlegendentry{${\boldsymbol x}^\star_{\rm W}$}

\addplot [smooth, mark=diamond*, Dandelion, mark size=4pt]
  table[row sep=crcr]{%
5	0.0304\\
10	0.1009\\
15	0.2217\\
20	0.4306\\
25	0.6698\\
30	0.9598\\
35	1.31\\
40	1.887\\
};
\addlegendentry{${ \boldsymbol x}^\star_{\rm LN}$}

\end{axis}
\end{tikzpicture}%
        \caption{Optimized interarrival times}
        \label{fig:subfig2a}
    \end{subfigure}
    \caption{Execution times for $\mbox{\sc scv}= 0.4$, time is in seconds.}
    \label{fig:main1}
\end{figure}

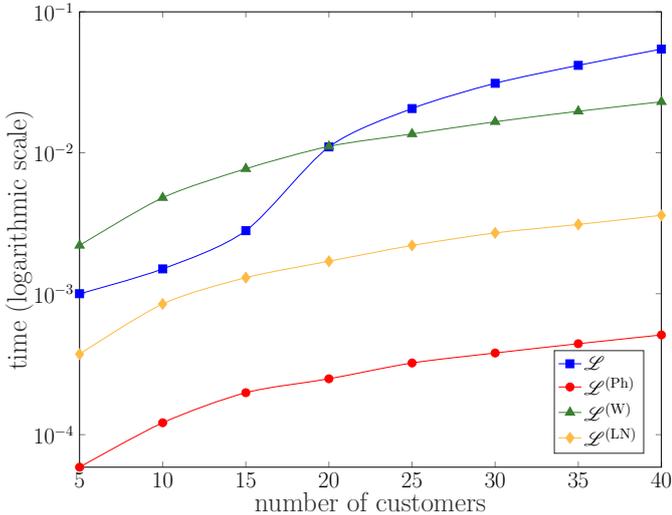
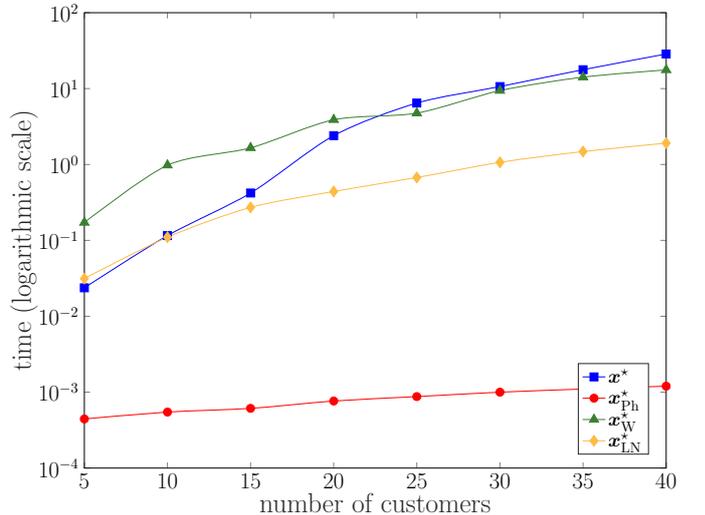
\begin{figure}[htbp]
    \begin{subfigure}{0.5\textwidth} 
        \centering
       \definecolor{mycolor1}{rgb}{1.00000,0.00000,1.00000}%
\begin{tikzpicture}[scale=0.5]

\begin{axis}[%
width=6.028in,
height=4.754in,
at={(1.011in,0.642in)},
scale only axis,
xmin=5,
xmax=40,
xtick={0,5,10,15,20,25,30,35,40},
xticklabel style={font=\fontsize{16}{14}\selectfont},
xlabel style={font=\color{white!15!black}},
xlabel={\fontsize{22}{14}\selectfont number of customers},
ymode=log,
ymin=5.8776e-05,
ymax=0.1,
yminorticks=false,
yticklabel style={font=\fontsize{16}{14}\selectfont},
ylabel style={font=\color{white!15!black}, yshift=10pt},
ylabel={\fontsize{22}{14}\selectfont time (logarithmic scale)},
legend style={at={(0.97,0.03)}, anchor=south east, legend cell align=left, align=left, draw=white!15!black, font=\fontsize{15}{14}\selectfont}
]
\addplot [smooth, mark=square*, blue, mark size=3pt]
  table[row sep=crcr]{%
5	0.001\\
10	0.0015\\
15	0.0028\\
20	0.011\\
25	0.0206\\
30	0.0311\\
35	0.0417\\
40	0.0544\\
};
\addlegendentry{$\mathscr{L}$}

\addplot [smooth,mark=*,red, mark size=3pt]
  table[row sep=crcr]{%
5	5.8776e-05\\
10	0.00012145\\
15	0.00019854\\
20	0.00024915\\
25	0.00032258\\
30	0.00037965\\
35	0.00044185\\
40	0.00050984\\
};
\addlegendentry{$\mathscr{L}^{(\rm Ph)}$}

\addplot [smooth, mark=triangle*, OliveGreen, mark size=4pt]
  table[row sep=crcr]{%
5	0.0022\\
10	0.0048\\
15	0.0077\\
20	0.0111\\
25	0.0136\\
30	0.0166\\
35	0.0197\\
40	0.023\\
};
\addlegendentry{$\mathscr{L}^{(\rm W)}$}

\addplot [smooth, mark=diamond*, Dandelion, mark size=4pt]
  table[row sep=crcr]{%
5	0.00037257\\
10	0.00084652\\
15	0.0013\\
20	0.0017\\
25	0.0022\\
30	0.0027\\
35	0.0031\\
40	0.0036\\
};
\addlegendentry{$\mathscr{L}^{(\rm LN)}$}

\end{axis}

\begin{axis}[%
width=7.778in,
height=5.833in,
at={(0in,0in)},
scale only axis,
xmin=0,
xmax=1,
ymin=0,
ymax=1,
axis line style={draw=none},
ticks=none,
axis x line*=bottom,
axis y line*=left
]
\end{axis}
\end{tikzpicture}%
        \caption{Loss functions}
        \label{fig:subfig1b}
    \end{subfigure}%
    \begin{subfigure}{0.5\textwidth} 
        \centering
       \definecolor{mycolor1}{rgb}{1.00000,0.00000,1.00000}%
\begin{tikzpicture}[scale=0.5]

\begin{axis}[%
width=6.028in,
height=4.754in,
at={(1.011in,0.642in)},
scale only axis,
xmin=5,
xmax=40,
xtick={0,5,10,15,20,25,30,35,40},
xticklabel style={font=\fontsize{16}{14}\selectfont},
xlabel style={font=\color{white!15!black}},
xlabel={\fontsize{22}{14}\selectfont number of customers},
ymode=log,
ymin=1e-04,
ymax=100,
yminorticks=false,
yticklabel style={font=\fontsize{16}{14}\selectfont},
ylabel style={font=\color{white!15!black}, yshift=10pt},
ylabel={\fontsize{22}{14}\selectfont time (logarithmic scale)},
legend style={at={(0.97,0.03)}, anchor=south east, legend cell align=left, align=left, draw=white!15!black, font=\fontsize{15}{14}\selectfont}
]
\addplot [smooth, mark=square*, blue, mark size=3pt]
  table[row sep=crcr]{%
5	0.0237\\
10	0.1156\\
15	0.4207\\
20	2.4006\\
25	6.4486\\
30	10.6372\\
35	17.746\\
40	28.5703\\
};
\addlegendentry{${\boldsymbol x}^\star$}

\addplot [smooth,mark=*,red, mark size=3pt]
  table[row sep=crcr]{%
5	0.00044368\\
10	0.00054671\\
15	0.00061211\\
20	0.00076458\\
25	0.00087434\\
30	0.001\\
35	0.0011\\
40	0.0012\\
};
\addlegendentry{${\boldsymbol x}^\star_{\rm Ph}$}

\addplot [smooth, mark=triangle*, OliveGreen, mark size=4pt]
  table[row sep=crcr]{%
5	0.1719\\
10	0.9784\\
15	1.6572\\
20	3.8857\\
25	4.7512\\
30	9.411\\
35	14.1549\\
40	17.6238\\
};
\addlegendentry{${\boldsymbol x}^\star_{\rm W}$}

\addplot [smooth, mark=diamond*, Dandelion, mark size=4pt]
  table[row sep=crcr]{%
5	0.0314\\
10	0.1102\\
15	0.272\\
20	0.441\\
25	0.6737\\
30	1.0675\\
35	1.481\\
40	1.9104\\
};
\addlegendentry{${\boldsymbol x}^\star_{\rm LN}$}

\end{axis}
\end{tikzpicture}%
        \caption{Optimized interarrival times}
        \label{fig:subfig2b}
    \end{subfigure}
    \caption{Execution times for $\mbox{\sc scv}= 0.7$, time is in seconds.}
    \label{fig:main2}
\end{figure}

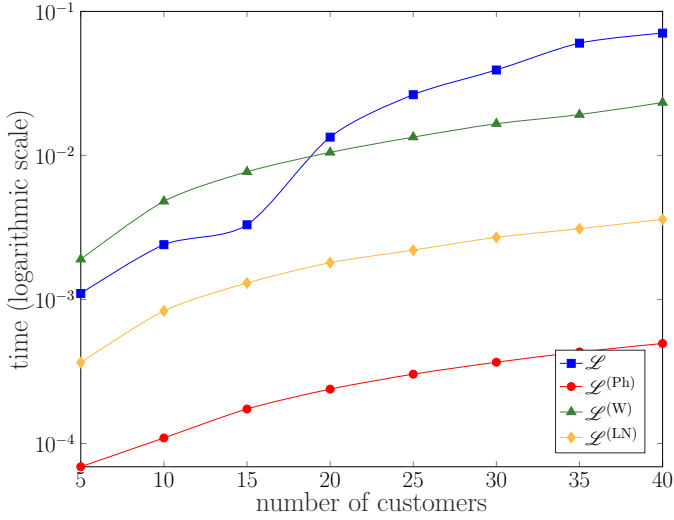
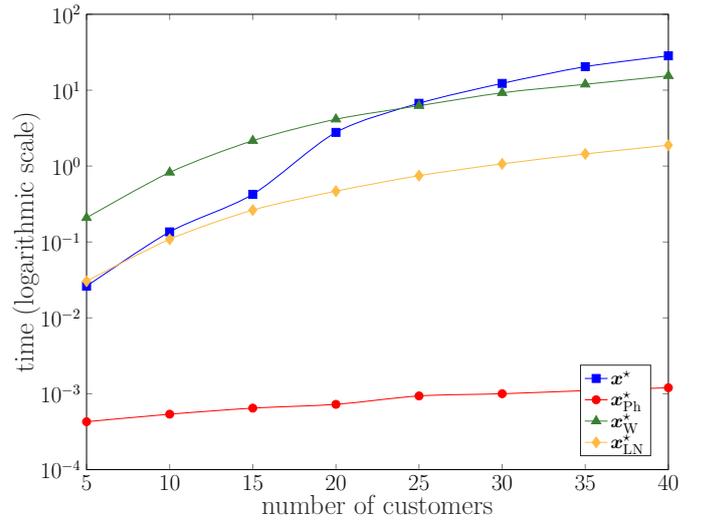
\begin{figure}[htbp]
    \begin{subfigure}{0.5\textwidth} 
        \centering
       \definecolor{mycolor1}{rgb}{1.00000,0.00000,1.00000}%
\begin{tikzpicture}[scale=0.5]

\begin{axis}[%
width=6.028in,
height=4.754in,
at={(1.011in,0.642in)},
scale only axis,
xmin=5,
xmax=40,
xtick={0,5,10,15,20,25,30,35,40},
xticklabel style={font=\fontsize{16}{14}\selectfont},
xlabel style={font=\color{white!15!black}},
xlabel={\fontsize{22}{14}\selectfont number of customers},
ymode=log,
ymin=6.8888e-05,
ymax=0.1,
yminorticks=false,
yticklabel style={font=\fontsize{16}{14}\selectfont},
ylabel style={font=\color{white!15!black}, yshift=10pt},
ylabel={\fontsize{22}{14}\selectfont time (logarithmic scale)},
legend style={at={(0.97,0.03)}, anchor=south east, legend cell align=left, align=left, draw=white!15!black, font=\fontsize{15}{14}\selectfont}
]
\addplot [smooth, mark=square*, blue, mark size=3pt]
  table[row sep=crcr]{%
5	0.0011\\
10	0.0024\\
15	0.0033\\
20	0.0134\\
25	0.0264\\
30	0.0392\\
35	0.0602\\
40	0.0707\\
};
\addlegendentry{$\mathscr{L}$}

\addplot [smooth,mark=*,red, mark size=3pt]
  table[row sep=crcr]{%
5	6.8888e-05\\
10	0.00010925\\
15	0.00017345\\
20	0.00023837\\
25	0.00030275\\
30	0.00036655\\
35	0.0004319\\
40	0.00049448\\
};
\addlegendentry{$\mathscr{L}^{(\rm Ph)}$}

\addplot [smooth, mark=triangle*, OliveGreen , mark size=4pt]
  table[row sep=crcr]{%
5	0.0019\\
10	0.0048\\
15	0.0077\\
20	0.0105\\
25	0.0134\\
30	0.0166\\
35	0.0192\\
40	0.0233\\
};
\addlegendentry{$\mathscr{L}^{(\rm W)}$}

\addplot [smooth, mark=diamond*, Dandelion, mark size=4pt]
  table[row sep=crcr]{%
5	0.00036467\\
10	0.00083198\\
15	0.0013\\
20	0.0018\\
25	0.0022\\
30	0.0027\\
35	0.0031\\
40	0.0036\\
};
\addlegendentry{$\mathscr{L}^{(\rm LN)}$}

\end{axis}

\begin{axis}[%
width=7.778in,
height=5.833in,
at={(0in,0in)},
scale only axis,
xmin=0,
xmax=1,
ymin=0,
ymax=1,
axis line style={draw=none},
ticks=none,
axis x line*=bottom,
axis y line*=left
]
\end{axis}
\end{tikzpicture}%
        \caption{Loss functions}
        \label{fig:subfig1c}
    \end{subfigure}%
    \begin{subfigure}{0.5\textwidth} 
        \centering
       \definecolor{mycolor1}{rgb}{1.00000,0.00000,1.00000}%
\begin{tikzpicture}[scale=0.5]

\begin{axis}[%
width=6.028in,
height=4.754in,
at={(1.011in,0.642in)},
scale only axis,
xmin=5,
xmax=40,
xtick={0,5,10,15,20,25,30,35,40},
xticklabel style={font=\fontsize{16}{14}\selectfont},
xlabel style={font=\color{white!15!black}},
xlabel={\fontsize{22}{14}\selectfont number of customers},
ymode=log,
ymin=0.0001,
ymax=100,
yminorticks=false,
yticklabel style={font=\fontsize{16}{14}\selectfont},
ylabel style={font=\color{white!15!black}, yshift=10pt},
ylabel={\fontsize{22}{14}\selectfont time (logarithmic scale)},
legend style={at={(0.97,0.03)}, anchor=south east, legend cell align=left, align=left, draw=white!15!black, font=\fontsize{15}{14}\selectfont}
]
\addplot [smooth, mark=square*, blue, mark size=3pt]
  table[row sep=crcr]{%
5	0.0261\\
10	0.1358\\
15	0.4244\\
20	2.7786\\
25	6.719\\
30	12.291\\
35	20.4043\\
40	28.3632\\
};
\addlegendentry{${\boldsymbol x}^\star$}

\addplot [smooth,mark=*,red, mark size=3pt]
  table[row sep=crcr]{%
5	0.00042792\\
10	0.00053788\\
15	0.00064661\\
20	0.00072551\\
25	0.00093613\\
30	0.001\\
35	0.0011\\
40	0.0012\\
};
\addlegendentry{${\boldsymbol x}^\star_{\rm Ph}$}

\addplot [smooth, mark=triangle*, OliveGreen, mark size=4pt]
  table[row sep=crcr]{%
5	0.2081\\
10	0.825\\
15	2.158\\
20	4.1389\\
25	6.2367\\
30	9.2407\\
35	11.9491\\
40	15.4659\\
};
\addlegendentry{${\boldsymbol x}^\star_{\rm W}$}

\addplot [smooth, mark=diamond*, Dandelion, mark size=4pt]
  table[row sep=crcr]{%
5	0.0304\\
10	0.1085\\
15	0.2632\\
20	0.4671\\
25	0.7479\\
30	1.0697\\
35	1.4453\\
40	1.8831\\
};
\addlegendentry{${\boldsymbol x}^\star_{\rm LN}$}

\end{axis}
\end{tikzpicture}%
        \caption{Optimized interarrival times}
        \label{fig:subfig2c}
    \end{subfigure}
    \caption{Execution times for $\mbox{\sc scv}= 1$, time is in seconds.}
    \label{fig:main3}
\end{figure}

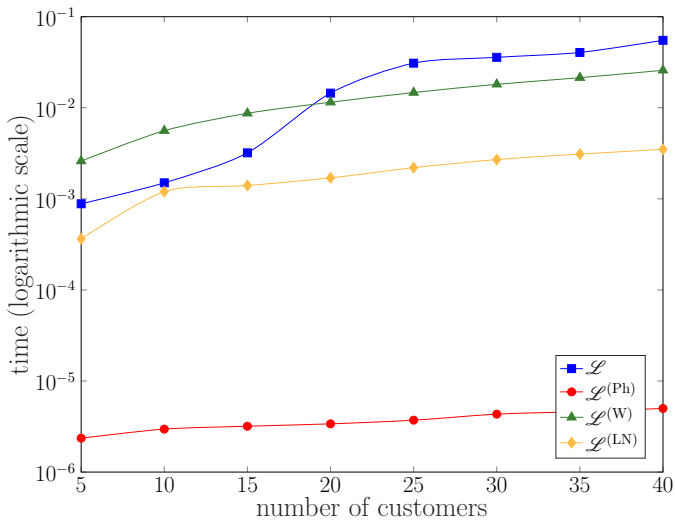
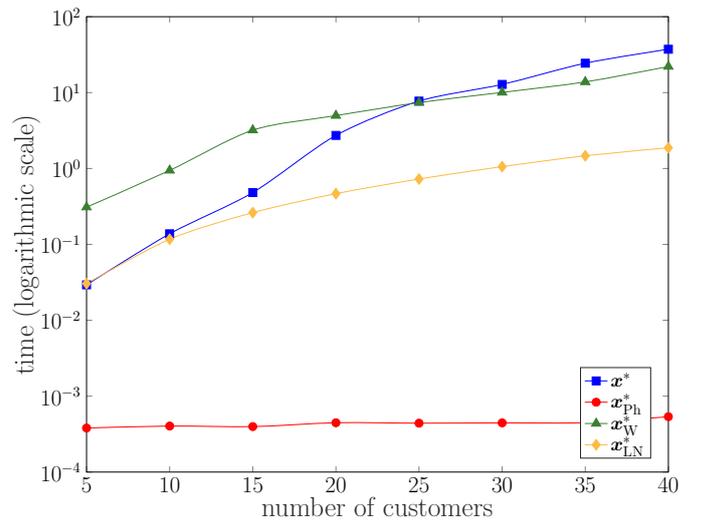
\begin{figure}[htbp]
    \begin{subfigure}{0.5\textwidth} 
        \centering
       \definecolor{mycolor1}{rgb}{1.00000,0.00000,1.00000}%
\begin{tikzpicture}[scale=0.5]

\begin{axis}[%
width=6.028in,
height=4.754in,
at={(1.011in,0.642in)},
scale only axis,
xmin=5,
xmax=40,
xtick={0,5,10,15,20,25,30,35,40},
xticklabel style={font=\fontsize{16}{14}\selectfont},
xlabel style={font=\color{white!15!black}},
xlabel={\fontsize{22}{14}\selectfont number of customers},
ymode=log,
ymin=1e-06,
ymax=0.1,
yminorticks=false,
yticklabel style={font=\fontsize{16}{14}\selectfont},
ylabel style={font=\color{white!15!black}, yshift=10pt},
ylabel={\fontsize{22}{14}\selectfont time (logarithmic scale)},
legend style={at={(0.97,0.03)}, anchor=south east, legend cell align=left, align=left, draw=white!15!black, font=\fontsize{15}{14}\selectfont}
]
\addplot [smooth, mark=square*, blue, mark size=3pt]
  table[row sep=crcr]{%
5	0.00088299\\
10	0.0015\\
15	0.0032\\
20	0.0145\\
25	0.0309\\
30	0.0358\\
35	0.0405\\
40	0.0551\\
};
\addlegendentry{$\mathscr{L}$}

\addplot [smooth,mark=*,red, mark size=3pt]
  table[row sep=crcr]{%
5	2.3547e-06\\
10	2.9641e-06\\
15	3.187e-06\\
20	3.3842e-06\\
25	3.7135e-06\\
30	4.3253e-06\\
35	4.6143e-06\\
40	4.9959e-06\\
};
\addlegendentry{$\mathscr{L}^{(\rm Ph)}$}

\addplot [smooth, mark=triangle*, OliveGreen, mark size=4pt]
  table[row sep=crcr]{%
5	0.0026\\
10	0.0056\\
15	0.0087\\
20	0.0115\\
25	0.0147\\
30	0.0181\\
35	0.0214\\
40	0.0258\\
};
\addlegendentry{$\mathscr{L}^{(\rm W)}$}

\addplot [smooth, mark=diamond*, Dandelion, mark size=4pt]
  table[row sep=crcr]{%
5	0.00036341\\
10	0.0012\\
15	0.0014\\
20	0.0017\\
25	0.0022\\
30	0.0027\\
35	0.0031\\
40	0.0035\\
};
\addlegendentry{$\mathscr{L}^{(\rm LN)}$}

\end{axis}
\end{tikzpicture}%
        \caption{Loss functions}
        \label{fig:subfig1d}
    \end{subfigure}%
    \begin{subfigure}{0.5\textwidth} 
        \centering
       \definecolor{mycolor1}{rgb}{1.00000,0.00000,1.00000}%
\begin{tikzpicture}[scale=0.5]

\begin{axis}[%
width=6.028in,
height=4.754in,
at={(1.011in,0.642in)},
scale only axis,
xmin=5,
xmax=40,
xtick={0,5,10,15,20,25,30,35,40},
xticklabel style={font=\fontsize{16}{14}\selectfont},
xlabel style={font=\color{white!15!black}},
xlabel={\fontsize{22}{14}\selectfont number of customers},
ymode=log,
ymin=1e-04,
ymax=100,
yminorticks=false,
yticklabel style={font=\fontsize{16}{14}\selectfont},
ylabel style={font=\color{white!15!black}, yshift=10pt},
ylabel={\fontsize{22}{14}\selectfont time (logarithmic scale)},
legend style={at={(0.97,0.03)}, anchor=south east, legend cell align=left, align=left, draw=white!15!black, font=\fontsize{15}{14}\selectfont}
]
\addplot [smooth, mark=square*, blue, mark size=3pt]
  table[row sep=crcr]{%
5	0.0293\\
10	0.1379\\
15	0.4821\\
20	2.7294\\
25	7.7266\\
30	12.8194\\
35	24.3923\\
40	37.2561\\
};
\addlegendentry{${\boldsymbol x}^*$}

\addplot [smooth,mark=*,red, mark size=3pt]
  table[row sep=crcr]{%
5	0.00038048\\
10	0.00040399\\
15	0.00039733\\
20	0.00044655\\
25	0.00044074\\
30	0.00044495\\
35	0.00044828\\
40	0.0005367\\
};
\addlegendentry{${\boldsymbol x}^*_{\rm Ph}$}

\addplot [smooth, mark=triangle*, OliveGreen, mark size=4pt]
  table[row sep=crcr]{%
5	0.3094\\
10	0.9438\\
15	3.2037\\
20	4.9744\\
25	7.3962\\
30	10.0555\\
35	13.8248\\
40	22.0529\\
};
\addlegendentry{${\boldsymbol x}^*_{\rm W}$}

\addplot [smooth, mark=diamond*, Dandelion, mark size=4pt]
  table[row sep=crcr]{%
5	0.0304\\
10	0.117\\
15	0.262\\
20	0.4657\\
25	0.7267\\
30	1.0569\\
35	1.4702\\
40	1.8729\\
};
\addlegendentry{${\boldsymbol x}^*_{\rm LN}$}

\end{axis}
\end{tikzpicture}%
        \caption{Optimized interarrival times}
        \label{fig:subfig2d}
    \end{subfigure}
    \caption{Execution times for $\mbox{\sc scv}= 1.3$, time is in seconds.}
    \label{fig:main4}
\end{figure}

After having assessed the accuracy of our approximations, in this section we analyze the computation time.
We do so by considering $n \in\{5,10,\ldots,40\}$ clients, for $\sigma^2:=\sigma_1^2=\cdots=\sigma_n^2 \in\{0.4,0.7,1.0,1.3\}$ and $\beta_1=\cdots=\beta_n=1$. 
We measure the computation time of
 \begin{itemize}
     \item[(A)]~The loss functions ${\mathscr L}({\boldsymbol x}_{\rm Eq}(y))$, ${\mathscr L}^{(\rm Ph)}({\boldsymbol x}_{\rm Eq}(y))$, ${\mathscr L}^{(\rm W)}({\boldsymbol x}_{\rm Eq}(y))$, and ${\mathscr L}^{(\rm LN)}({\boldsymbol x}_{\rm Eq}(y))$ for $y=1.5$.
    \item[(B)]~The optimized interarrival times ${\boldsymbol x}^\star$, ${\boldsymbol x}_{\rm Ph}^\star$, ${\boldsymbol x}_{\rm W}^\star$, and ${\boldsymbol x}_{\rm LN}^\star$.
\end{itemize}
Figures \ref{fig:main1}--\ref{fig:main4} display the corresponding execution times, each for one of the four values of the {\sc scv}. 
These figures clearly show that, in terms of execution time, our phase-type approximation outperforms the alternatives by a huge margin. 
This is true not just for the evaluation of the loss function ${\mathscr L}^{\rm (Ph)}({ x}_{\rm Eq})$, but also for computing the optimized interarrival times ${\boldsymbol x}^\star_{\rm Ph}$.
 
 We proceed by having a closer look at the computation times displayed by Figures \ref{fig:main1}--\ref{fig:main4}.
 We first consider the execution times required for the evaluation of the loss function (left panels in Figures \ref{fig:main1}--\ref{fig:main4}). 
 While, due to the matrix operations involved, the evaluation of ${\mathscr L}({\boldsymbol x}_{\rm Eq}(y))$ scales essentially cubically in the number of customers $n$, the computation time of ${\mathscr L}^{(\rm Ph)}({\boldsymbol x}_{\rm Eq}(y))$ is only {\it linear} in the number of customers $n$. 
 In this instance this yields a speedup of more than a factor $10^3$; in all other instances considered we observed gains of the same order. 
 The other two approximations, using the Weibull and Lognormal fit, also scale linearly in $n$ (as witnessed by the parallel curves in the left panels of Figures \ref{fig:main1}--\ref{fig:main4}), but with substantially larger proportionality constants --- this is a consequence of the necessity to evaluate certain non-explicit functions when performing the two-moments fit (namely the incomplete Gamma function and the standard Normal complementary distribution function), and in addition in the Weibull case a bisection is needed to identify the $\alpha$-parameter. 
 These conclusions are consistent across the various values of {\sc scv}. 

 Essentially the same conclusions carry over to the execution times for the optimized interarrival times  (right panels in Figures \ref{fig:main1}--\ref{fig:main4}). In general, ${\boldsymbol x}^\star_{\rm Ph}$ performs (by far) best, with an execution time that is essentially linear in $n$.  
 In this example, the computation times for ${\boldsymbol x}^\star_{\rm Ph}$ are even sublinear, but that is a consequence of the fact that we have worked with homogeneous customers. 
 Such homogeneous scenarios are of lower complexity, because the resulting optimized schedule is fairly equidistant, except for the interarrival times of the first and last few customers (in the literature known as the `dome shape').
 
 The figures reveal that when using the phase-type fit, the evaluation of the loss function for 40 customers takes just between $10^{-4}$ and $10^{-3}$ s. 
 Computing the corresponding optimal schedule, which minimizes the loss function over the vector of arrival times, takes roughly $10^{-3}$ s. 
 These low execution times are  particularly promising in the context of problems in which one not only aims at determining optimal arrival times, but in addition wants to identify the customers' optimal order; cf.\ the sequencing problem analyzed in \cite{de2021} and the routing problem analyzed in \cite{bekker2023}.

 Additional experiments have shown that the above conclusions remain valid for heterogeneous scenarios. 
 In the experiments underlying Table~\ref{heterogeneous}, the execution time when using the phase-type evaluation ${\mathscr L}$ was $1.5\cdot 10^{-1}$ sec, whereas the approximation based on the phase-type fit took $3.0\cdot 10^{-4}$ sec, the approximation based on the Weibull fit took $2.3\cdot 10^{-2}$~sec, and the approximation based on the Lognormal fit took $3.8\cdot 10^{-3}$~sec. 

 \subsection{Conclusions from numerical experiments}
The three two-moments fits that have been used in our approximations are roughly equally accurate, except in certain specific regimes. Across all instances considered, chosen such that they cover all relevant scenarios, the phase-type fit performs well (but typically excellently). 
As this phase-type fit by far outperforms the other fits in terms of execution time, it is the approximation of our choice.

\section{Discussion and concluding remarks}\label{DIS}
We start this section by commenting on the application potential of our approach. 
In this respect, we mention that it facilitates an efficient and accurate approximation  of the expected waiting and idle times, given the arrival times of the individual customers, for generally (not necessarily homogeneous) service-time distributions.
This allows the evaluation of the loss function \eqref{OBJ0} that is generally used in the appointment scheduling literature. 

An exact algorithm that outputs \eqref{OBJ0} is only available for the case when the service times have phase-type distributions, but it tends to be time consuming.
In particular if the number of customers gets large, this approach becomes prohibitively slow, as was confirmed by the experiments discussed in Section~\ref{comp_speed}. 
We had been facing this problem in our predecessor project \cite{bekker2023}, in which we aimed to use appointment scheduling methodology in a routing context. The procedures developed in the present paper can be used to substantially speed up such routing algorithms.

We developed three approximation algorithms, in which the iterands in the Lindley recursion are approximated via two-moments fits (using a low-dimensional phase-type fit, a Weibull fit and a Lognormal fit). We thoroughly tested the approximations through a carefully designed series of experiments (a selection of which has been presented in Section \ref{NUM}), where we chose the parameters so as to maximally cover the parameter space. Based on these experiments, we recommend the use of the phase-type fit. While the other fits are roughly equally robust, in that they provide essentially equally accurate results (except in certain specific regimes), the phase-type fit is, by a huge margin, superior in terms of execution time. Using this phase-type fit, the evaluation of the loss function for 40 customers takes between $10^{-4}$ and $10^{-3}$ sec on a standard MacBook Pro. Computing the corresponding optimal schedule, where the loss function is optimized over the vector of arrival times, takes about as little as $10^{-3}$ sec.

\vb

There are various directions for future research. In our context of appointment schedules, the interarrival times are deterministic, but these can be taken stochastic as well. A second research theme concerns the extension to multi-server queues, where one should work with the Kiefer-Wolfowitz recursion rather than the Lindley recursion \cite[Chapter XII]{asmussen2003}. Third, one could think of network extensions, for instance a tandem network; this has the intrinsic complication that interdeparture times are not independent. A fourth research direction could be the case of single-server queues with a more involved service mechanism, such as priority queues or parallel queues.


\bibliographystyle{abbrv}
\bibliography{references}

\begin{thebibliography}{10}

\bibitem{abate1995}
J.~Abate and W.~Whitt.
\newblock Numerical inversion of {L}aplace transforms of probability
  distributions.
\newblock {\em ORSA Journal on Computing}, 7(1):36--43, 1995.

\bibitem{ahmadi2017}
A.~Ahmadi-Javid, Z.~Jalali, and K.~J. Klassen.
\newblock Outpatient appointment systems in healthcare: {A} review of
  optimization studies.
\newblock {\em European Journal of Operational Research}, 258(1):3--34, 2017.

\bibitem{asmussen2003}
S.~Asmussen.
\newblock {\em Applied Probability and Queues}.
\newblock Springer, New York, 2nd edition, 2003.

\bibitem{asmussen2007}
S.~Asmussen and P.~W. Glynn.
\newblock {\em Stochastic Simulation: Algorithms and Analysis}, volume~57.
\newblock Springer, 2007.

\bibitem{bekker2023}
R.~Bekker, Bharti, L.~Lan, and M.~Mandjes.
\newblock A queueing-based approach for integrated routing and appointment
  scheduling.
\newblock {\em Submitted}, 2023.

\bibitem{bux1979}
W.~Bux.
\newblock Single-server queues with general interarrival and phase-type service
  time distributions.
\newblock In {\em Proceedings of the 9th International Teletraffic Congress,
  Torremolinos}, 1979.

\bibitem{de2021}
M.~A. de~Kemp, M.~Mandjes, and N.~Olver.
\newblock Performance of the smallest-variance-first rule in appointment
  sequencing.
\newblock {\em Operations Research}, 69(6):1909--1935, 2021.

\bibitem{de1989}
A.~de~Kok.
\newblock A moment-iteration method for approximating the waiting-time
  characteristics of the {GI/G/1} queue.
\newblock {\em Probability in the Engineering and Informational Sciences},
  3(2):273--287, 1989.

\bibitem{ho1992}
C.-J. Ho and H.-S. Lau.
\newblock Minimizing total cost in scheduling outpatient appointments.
\newblock {\em Management Science}, 38(12):1750--1764, 1992.

\bibitem{klosterhalfen2018}
S.~T. Klosterhalfen, F.~Holzhauer, and M.~Fleischmann.
\newblock Control of a continuous production inventory system with production
  quantity restrictions.
\newblock {\em European Journal of Operational Research}, 268(2):569--581,
  2018.

\bibitem{kuiper2015}
A.~Kuiper, B.~Kemper, and M.~Mandjes.
\newblock A computational approach to optimized appointment scheduling.
\newblock {\em Queueing Systems}, 79:5--36, 2015.

\bibitem{kuiper2017}
A.~Kuiper, M.~Mandjes, and J.~de~Mast.
\newblock Optimal stationary appointment schedules.
\newblock {\em Operations Research Letters}, 45(6):549--555, 2017.

\bibitem{kuiper2023}
A.~Kuiper, M.~Mandjes, J.~de~Mast, and R.~Brokkelkamp.
\newblock A flexible and optimal approach for appointment scheduling in
  healthcare.
\newblock {\em Decision Sciences}, 54(1):85--100, 2023.

\bibitem{lewis2008}
A.~L. Lewis and E.~Mordecki.
\newblock Wiener-{H}opf factorization for {L}{\'e}vy processes having positive
  jumps with rational transforms.
\newblock {\em Journal of Applied Probability}, 45(1):118--134, 2008.

\bibitem{neuts1981}
M.~F. Neuts.
\newblock {\em Matrix-Geometric Solutions in Stochastic Models: An Algorithmic
  Approach}.
\newblock The Johns Hopkins Press, Baltimore, 1981.

\bibitem{prabhu1998}
N.~U. Prabhu.
\newblock {\em Stochastic Storage Processes: Queues, Insurance Risk, and Dams,
  and Data Communication}.
\newblock Number~15. Springer Science \& Business Media, 1998.

\bibitem{robinson2003}
L.~W. Robinson and R.~R. Chen.
\newblock Scheduling doctors' appointments: {O}ptimal and empirically-based
  heuristic policies.
\newblock {\em IIE Transactions}, 35(3):295--307, 2003.

\bibitem{seelen1984}
L.~Seelen and H.~C. Tijms.
\newblock Approximations for the conditional waiting times in the {GI/G/c}
  queue.
\newblock {\em Operations Research Letters}, 3(4):183--190, 1984.

\bibitem{tijms1986}
H.~C. Tijms.
\newblock {\em Stochastic Modelling and Analysis: A Computational Approach}.
\newblock John Wiley \& Sons, Inc., Chichester \& New York, 1986.

\bibitem{wagner2004}
M.~Wagner and S.~R. Smits.
\newblock A local search algorithm for the optimization of the stochastic
  economic lot scheduling problem.
\newblock {\em International Journal of Production Economics}, 90(3):391--402,
  2004.

\bibitem{wang1997}
P.~P. Wang.
\newblock Optimally scheduling $n$ customer arrival times for a single-server
  system.
\newblock {\em Computers \& Operations Research}, 24(8):703--716, 1997.

\bibitem{zhan2021}
Y.~Zhan, Z.~Wang, and G.~Wan.
\newblock Home service routing and appointment scheduling with stochastic
  service times.
\newblock {\em European Journal of Operational Research}, 288(1):98--110, 2021.

\bibitem{zijm1994}
H.~Zijm and G.-J. van Houtum.
\newblock On multi-stage production/inventory systems under stochastic demand.
\newblock {\em International Journal of Production Economics},
  35(1-3):391--400, 1994.

\end{thebibliography}

\appendix
\section{Excess distribution for mixed Erlang}\label{APP}
This appendix presents some results on the excess $X_x:=(X-x)^+$ when $X\sim {\rm ME}(p,k,\mu).$
As mentioned in the body of the paper, ${\mathbb E}\,X_x$ can be evaluated through a direct computation, but this appendix gives an elegant probabilistic derivation.

First, consider the case that $X$ follows an exponential distribution with rate $\mu$, providing the essence of the approach. Due to the lack of memory of $X$, we have \[\ee \,X_x = \pp(X>x)\, \ee\left[ X-x \,|\, X>x\right] = e^{-\mu x} \frac{1}{\mu}.\]
Similarly, for the second moment, we have
\[\ee (X_x)^2 = \pp(X>x) \,\ee\left[ (X-x)^2 \,|\, X>x\right] = e^{-\mu x} \frac{2}{\mu^2}.\] If $X\sim {\rm ME}(p,k,\mu)$, then observe that the number of exponential phases that are completed during $(0,x)$ follows a Poisson process. Moreover, if $i =0,\ldots,k-1$ exponential phases have been completed during $(0,x)$, then $X-x$ follows the mixed Erlang distribution $X\sim {\rm ME}(p,k-i,\mu)$.
Hence, combining the above and conditioning on the number of completed exponential phases, yields
\begin{align*}
\ee \,X_x & = \sum_{i=0}^{k-1} e^{-\mu x} \frac{(\mu x)^i}{i!} \left( (1-p) \; \frac{k-i}{\mu} + p \; \frac{k-1-i}{\mu} \right) \\
 & = \sum_{i=0}^{k-1} e^{-\mu x} \frac{(\mu x)^i}{i!} \frac{k-p-i}{\mu} = \frac{1}{\mu} e^{-\mu x} \left((k-p-\mu x) \sum_{i=0}^{k-2} \frac{(\mu x)^i}{i!} + (k-p) \frac{(\mu x)^{k-1}}{(k-1)!} \right),
\end{align*} 
where the third equality follows from splitting the final term in the second equation and some elementary rewriting. 
Alternatively, the expected excess may be expressed in terms of the incomplete Gamma integral using the well-known identity 
\[\frac{\Gamma(k,t)}{(k-1)!} = e^{-t} \sum_{i=0}^{k-1} \frac{t^i}{i!}.\] Specifically, we may then write
\begin{equation*}
\ee \,X_x = \frac{ (k-p-\mu x)}{\mu (k-2)!} \Gamma(k-1, \mu x) + \frac{(k-p)}{\mu (k-1)!} (\mu x)^{k-1} e^{-\mu x}.  
\end{equation*}
Observe that this is consistent with the definition of $M(k,x,\mu)$, whereas the expression for $r_{i+1}$ in \eqref{RiME} can be slightly simplified. 

Finally, for the second moment of $R$ we may apply similar arguments. Recall that $\ee X^2 = k(k+1-2p)/\mu^2$. Conditioning on the number of completed exponential phases during $(0,x)$ again, we obtain
\begin{align*}
\ee (X_x)^2 
 & = \sum_{i=0}^{k-1} e^{-\mu x} \frac{(\mu x)^i}{i!} \frac{(k-i)(k-2p+1-i)}{\mu^2}  = \frac{1}{\mu^2} e^{-\mu x}  \left[ c_1 \sum_{i=0}^{k-1} \frac{(\mu x)^i}{i!} + c_2 \frac{(\mu x)^{k-1}}{(k-1)!} \right],
\end{align*} 
where $c_1:=(k-\mu x)^2+k-2p(k+\mu x)$ and $c_2:=\mu x (1+k+2p-\mu x)$. Similar to what we have done for the first moment, we may write
\begin{equation*}
\ee (X_x)^2 = \frac{c_1}{\mu^2 (k-1)!} \Gamma(k, \mu x) + \frac{c_2}{\mu^2 (k-1)!} (\mu x)^{k-1} e^{-\mu x}.  
\end{equation*}
This expression is consistent with $S(k,x,\mu)$, whereas the expression for $v_{i+1}$ in \eqref{ViME} can be slightly simplified.

\section{Two-moment fit}\label{Explicit_FIT}
 
In this appendix, we provide the two-moment fit to the mean $r_i$ and variance $v_i$ for hyperexponential, mixed Erlang, Weibull, and Lognormal distribution functions. 
In this appendix we denote $\mbox{\sc scv}_i=v_i/r_i^2$ for the squared coefficient of variation pertaining to the distribution of client $i$'s service time.

\subsubsection*{Hyperexponential fit}
In case $R_i \sim  {\rm HE}(\alpha,\mu_1,\mu_2)$, corresponding to the situation that $\mbox{\sc scv}_i \geqslant 1$, we need to find three parameters with only two constraints, and thus we have one degree of freedom left. 
To obtain a unique parametrization, we impose the extra condition of `balanced means' from \cite{tijms1986}, i.e., we require that $\mu_1=2\alpha \mu$ and $\mu_2=2(1-\alpha)\mu$ for some $\mu>0$. 
It can be verified (see \cite{tijms1986}) that one should take
\begin{equation} \label{eqn:H2-fit} 
\alpha = \frac{1}{2}
  \left(1 + \sqrt{{\frac{\mbox{\scv}_i-1}{\mbox{\scv}_i+1}}}\right),
 \quad \mu_1 = \frac{2\alpha}{r_i}, 
 \quad \mu_2 = \frac{2(1-\alpha)}{r_i}.
\end{equation}

\subsubsection*{Mixed-Erlang fit}
In case $R_i \sim {\rm ME}(p,k,\mu)$, corresponding to the situation that $\mbox{\sc scv}_i < 1$, we first select the appropriate $k$ such that $\mbox{\scv}_i\in \left(k^{-1},(k-1)^{-1}\right]$. Subsequently, we take 
\begin{equation} \label{eqn:ME-fit} 
p = \frac{1}{1+\mbox{\scv}_i} \left( k \cdot \mbox{\scv}_i - \sqrt{k (1+\mbox{\scv}_i) - k^2 \cdot\mbox{\scv}_i} \right)
    \quad {\rm and} \quad
    \mu = \frac{k-p}{r_i}.
\end{equation}

\subsubsection*{Weibull fit}\label{FIT_Weib}
For $R_i \sim {\rm W}(\lambda,\alpha)$, we find $\alpha\equiv \alpha(r_i,v_i)$ from the equation
\[\frac{\Gamma(1+2/\alpha)-(\Gamma(1+1/\alpha))^2}{(\Gamma(1+1/\alpha))^2}= \mbox{\scv}_i,\]
and then 
\[\lambda=r_i^{-1}\cdot\Gamma\Big(1+\frac{1}{\alpha(r_i,v_i)}\Big) .\]

\subsubsection*{Lognormal fit}\label{FIT_Logn}
For $R_i \sim {\rm LN}(\mu,\tau^2)$, it follows after an elementary calculation that the two-moment fit equals
\begin{equation} \label{eqn:LN-fit}
    \mu=\log{\Big(\frac{r_i^2}{\sqrt{v_i+r_i^2}}\Big)}~~~\text{and}~~~\tau^2=\log{\Big(1+\frac{v_i}{r_i^2}\Big)}. 
\end{equation}

\end{document}